\documentclass{article}

\usepackage{silence}
\usepackage[T1]{fontenc}
\usepackage[utf8]{inputenc}

\usepackage[margin=1in]{geometry}
\usepackage{amssymb, amsmath, amsthm}

\usepackage{tikz}
\usetikzlibrary{calc,matrix,arrows,trees,positioning}

\newif\ifpdf \pdffalse \ifx\pdfoutput\undefined\else\ifx\pdfoutput\relax\else\ifnum\pdfoutput<1 \else\pdftrue\fi\fi\fi
\ifpdf
\usepackage[pdfencoding=auto,psdextra,backref=section,linktoc=all]{hyperref}
\else
\usepackage[hypertex,backref=section,linktoc=all]{hyperref}
\fi

\usepackage[capitalize,nosort]{cleveref}

\crefformat{equation}{(#2#1#3)}
\crefname{figure}{Figure}{Figures}

\numberwithin{equation}{section}

\theoremstyle{definition}
\newtheorem{theorem}[equation]{Theorem}
\newtheorem{corollary}[equation]{Corollary}
\newtheorem{lemma}[equation]{Lemma}
\newtheorem{notation}[equation]{Notation}
\newtheorem{proposition}[equation]{Proposition}

\newtheorem{definition}[equation]{Definition}
\newtheorem{example}[equation]{Example}

\newtheorem{remark}[equation]{Remark}

\newcommand{\PSh}{\mathsf{PSh}}
\newcommand{\Sh}{\mathsf{Sh}}

\newcommand{\sk}{{\rm sk}}

\newcommand{\R}{\mathbb{R}}
\newcommand{\NN}{\mathbf{N}} 

\newcommand{\op}{{\rm op}}

\newcommand{\Sp}{\mathsf{S}} 
\newcommand{\sSp}{\mathsf{sS}} 
\newcommand{\Sets}{\mathsf{Sets}}

\def\rdf{{\bf R}}
\newcommand{\Hom}{\mathop{\rm Hom}\nolimits}
\newcommand{\RHom}{\rdf\mathop{\rm Hom}\nolimits}
\newcommand{\sd}{\mathop{\rm sd}\nolimits}
\newcommand{\Sd}{\mathop{\rm Sd}\nolimits}
\newcommand{\Ex}{\mathop{\rm Ex}\nolimits}
\newcommand{\map}{\mathop{\rm map}\nolimits}
\newcommand{\Rmap}{\rdf\mathop{\rm map}\nolimits}

\newcommand{\hocolim}{\mathop{\rm hocolim}\limits}
\newcommand{\holim}{\mathop{\rm holim}\limits}
\newcommand{\Sing}{\mathop{\rm Sing}\limits}

\def\bD{\mathsf{D}}
\newcommand{\A}{\mathbb{A}}
\newcommand{\sC}{\mathbf{C}}
\newcommand{\sU}{\mathbf{U}}
\newcommand{\man}{\mathsf{Man}}
\newcommand{\sB}{\mathbf{B}}
\newcommand{\tB}{{\rm B}}
\newcommand{\sD}{\mathbf{D}}
\newcommand{\id}{{\rm id}}
\newcommand{\pre}{{\rm pre}}
\def\bDelta{{\bf\Delta}}
\def\inj{{\rm inj}}
\def\int{\mathop{\rm int}}
\def\cU{\mathcal{U}}
\def\cV{\mathcal{V}}

\mathcode`:="603A 
\mathchardef\colon="303A
\def\coloneq{\colon=}

\newcommand{\downmap}[2]{
\vcenter{\halign{
&\hfil$##$\hfil\cr
#1\cr\downarrow\cr#2\cr}}
}

\newcommand\leftleftleftarrows{
\mathrel{\substack{\textstyle\leftarrow\\[-0.5ex]
                      \textstyle\leftarrow \\[-0.5ex]
                      \textstyle\leftarrow}}
}

\begin{document}

\title{Classifying spaces of infinity-sheaves}
\author{Daniel Berwick-Evans\\
Department of Mathematics, University of Illinois Urbana--Champaign%
\thanks{
\url{https://danbe.web.illinois.edu/};
danbe@illinois.edu
}
\and
Pedro Boavida de Brito\\
Department of Mathematics, Instituto Superior Técnico, University of Lisbon%
\thanks{
\url{https://www.math.tecnico.ulisboa.pt/\~pbrito/};
pedrobbrito@tecnico.ulisboa.pt
}
\and
Dmitri Pavlov\\
Department of Mathematics and Statistics, Texas Tech University%
\thanks{
\url{https://dmitripavlov.org/}
}}

\date{}

\maketitle


\begin{abstract}
We prove that the set of concordance classes of sections of an $\infty$-sheaf on a manifold is representable, extending a theorem of Madsen and Weiss for sheaves of sets.
This is reminiscent of an $h$-principle in which the role of isotopy is played by concordance.
As an application, we offer an answer to the question: what does the classifying space of a Segal space classify?
\end{abstract}

\tableofcontents

\section{Introduction}

Let $F$ be an $\infty$-sheaf (alias homotopy sheaf, see \cref{defn:sheaf}) on $\man$, the site of finite-dimensional smooth manifolds without boundary and smooth maps.
For a manifold $M$, an element of $F(M \times \R)$ is called a \emph{concordance}.
Two elements $\sigma_0$ and $\sigma_1$ in $F(M)$ are said to be \emph{concordant} (and we write $\sigma_0\sim_c \sigma_1$) if there exists a concordance whose restriction to $M \times \{k\}$ is $\sigma_k$ for $k = 0,1$.

Concordance is an equivalence relation, and a familiar one in many situations.
Here are three examples.
When $F=C^\infty(-,N)$, maps are concordant if and only if they are smoothly homotopic.
For the sheaf of closed differential $n$-forms, two sections (i.e., closed $n$-forms) are concordant if and only if they are cohomologous.
For the stack of vector bundles, a pair of vector bundles are concordant if and only if they are isomorphic.
In these three cases, concordance classes have a well-known description in terms of homotopy classes of maps into a space,
namely the space underlying~$N$, the Eilenberg--MacLane space~${\rm K}(\R,n)$, and the space~${\rm BO}(n)$, respectively.
In this paper we generalize these classical representability results:
concordance classes of sections of \emph{any} $\infty$-sheaf $F$ is represented by a space $\tB F$, which we call the \emph{classifying space} of $F$.

We now assemble the ingredients to state our results precisely.
We denote by $\A^n$ the smooth extended simplex, that is, the subspace of $\R^{n+1}$ whose coordinates sum to one.
By varying $n$, this defines a cosimplicial object in $\man$.
Define a presheaf
$$\sB F(M) \coloneq |[k] \mapsto F(M \times \A^k)|$$
with values in spaces (i.e., simplicial sets), where $|{-}|$ denotes the homotopy colimit of the simplicial space.
The construction $\sB$ is a form of localization; it is the universal way to render the maps $F(M) \to F(M \times \A^1)$ invertible for all $M$ and $F$.
It is a familiar construction in the motivic literature,
for example in the work of Morel--Voevodsky~\cite{Morel.Voevodsky} (who call it $Sing$),
but it has also appeared in the context of geometric topology in the papers of
Waldhausen~\cite{Waldhausen}, Weiss--Williams~\cite{Weiss.Williams}, and Madsen--Weiss~\cite{Madsen.Weiss}.
The link between $\sB F$ and the concordance relation~$\sim_c$ is the bijection $\pi_0 \sB F(M)\cong \pi_0 F(M)/ \sim_c$.

\medskip
Define the \emph{classifying space} $\tB F$ as a Kan complex replacement of $\sB F(*)$ and denote by $\Sing M$ the usual singular simplicial set of $M$.
Our main result is:

\begin{theorem}\label{mainthm}
Let $F$ be an $\infty$-sheaf on $\man$.
There is an evaluation map
$$(\sB F)(M) \rightarrow \Rmap(\Sing M, \tB F)$$
which is a natural weak equivalence of spaces for every manifold $M$.
\end{theorem}

It is not difficult to show---essentially by a variant of Brown's representability theorem; see \cref{sec:concinv}---that \cref{mainthm} is equivalent to the following.
\begin{theorem}\label{sheafthm}
If $F$ is an $\infty$-sheaf, then $\sB F$ is an $\infty$-sheaf.
\end{theorem}

These statements may be regarded as analogues of the $h$-principle, where the usual relation of \emph{isotopy} is replaced by that of \emph{concordance}.
Here we have in mind the strand of the $h$-principle that gives conditions (e.g., microflexibility) which guarantee that an isotopy invariant functor (e.g., a sheaf) is an $\infty$-sheaf.
The relation of concordance is more severe than that of isotopy, and this explains why the hypotheses are less restrictive than those of typical $h$-principles, e.g., there are no dimension restrictions, open versus closed manifolds, etc.

\medskip
Just as with the $h$-principle, the key step in our proof involves verifying certain fibration properties.
As such,
a significant part of the paper is a study of weak lifting properties for maps of simplicial spaces.
We introduce the notion of \emph{weak Kan fibration} of simplicial spaces and simplicial sets.
A crucial result shows that weak Kan fibrations are realization fibrations (see \cref{defn:realfib} and \cref{thm:wkanmain}); this implies that geometric realization is stable under homotopy pullback along weak Kan fibrations.

\medskip
We emphasize that these results---and hence \cref{sheafthm}---do not follow from formal considerations.
There are simple counterexamples in the category of schemes as in the $\A^1$-homotopy theory of Morel--Voevodsky~\cite[\S 3, Example 2.7]{Morel.Voevodsky}.
The $\infty$-sheaf property is a homotopy limit condition whereas $\sB$ involves geometric realization, a homotopy colimit.
Commuting these is a subtle issue.
This is where we use the weak Kan property to prove \cref{sheafthm}.
The verification that $\sB F$ is weak Kan and certain restriction maps are weak Kan fibrations follows from a geometric argument about smooth maps (\cref{lem:relhorninc}).
We also rely on the existence of partitions of unity.

\medskip
Our main results improve on prior work of others, although the techniques we use differ.
The $\pi_0$-statement of \cref{mainthm} is due to Madsen--Weiss~\cite[Appendix~A]{Madsen.Weiss} when $F$ is a sheaf of \emph{sets} (or of discrete categories).
Although our theorem does extend the result of Madsen--Weiss from $\pi_0$ to $\pi_n$ in the case of sheaves of sets, the main objective of our work is to extend it from sheaves of sets to $\infty$-sheaves of spaces.
The argument of Madsen--Weiss shares some features with ours (in that certain locally constancy conditions along simplices of a triangulation are enforced), but does not extend to $\infty$-sheaves.
Moreover, unlike theirs, our arguments apply in the topological or PL category too: \cref{mainthm} remains true if we consider topological or PL manifolds instead of smooth manifolds.
In fact, our arguments simplify significantly in those cases (see \cref{sec:weakKanconc} for explanations).

Bunke--Nikolaus--V\"olkl~\cite[\S7]{Bunke.Nikolaus.Volkl} have proved a version of \cref{mainthm} for $\infty$-sheaves on compact manifolds with values in \emph{spectra}.
From the point of view of \cref{sheafthm}, this case is essentially formal since in a stable setting homotopy pullback squares are homotopy pushout squares, and so (for finite covers) the problem of commuting homotopy pullbacks with geometric realization disappears.
It has also been pointed out to us by John Francis that Ayala–Francis–Rozenblyum~\cite[2.3.16, 2.4.5]{Ayala.Francis.Rozenblyum} gave related results that are proved in the context of stratified spaces.
As we understand it, these results are both more general (they apply to \emph{stratified} spaces) and less general than ours and those of Madsen--Weiss (they apply to a certain class of \emph{isotopy} sheaves on groupoids).
This restricted class of sheaves is, from the point of view of \cref{mainthm}, too severe as it excludes many $\infty$-sheaves, even set-valued ones.

\medskip
Applications of \cref{mainthm} abound.
This stems from the fact that we not only prove abstract representability, but also give a formula for the representing space.
This formula can be identified with classical constructions.
Two illustrative examples, connecting back to the beginning of this introduction, are the classical de Rham theorem and the classification of vector bundles (with or without connection).
These are obtained by applying the main theorem to the sheaf of differential $n$-forms and the stack of vector bundles (sheaves of sets and stacks are examples of $\infty$-sheaves).
In \cref{sec:class}, we discuss a further application: a classification of $\sC$-bundles, where $\sC$ is a Segal space.
More recently, the third author proved \cite[Theorem~13.8]{Pavlov.Diffeo} a generalization of \cref{mainthm} for presheaves with values in model categories Quillen equivalent
to model categories of algebras over simplicial algebraic theories (such as chain complexes, connective spectra, and various flavors of connective ring spectra).

\medskip

We mention here another consequence of \cref{mainthm}.
Let $\bD$ denote the full subcategory of $\man$ spanned by $\R^n$ with $n \geq 0$.
The $\infty$-categories $\Sh(\man)$ and $\Sh(\bD)$ of $\infty$-sheaves on $\man$ and $\bD$, respectively (with respect to the usual open covers by codimension zero embeddings), are examples of $\infty$-toposes, as are the slice $\infty$-categories $\Sh(\man)/F$ for an $\infty$-sheaf $F$ on $\man$, and $\Sh(\bD)/F$ for an $\infty$-sheaf $F$ on $\bD$.

\begin{proposition}
The functor $\tB$ from $\Sh(\man)$ to spaces is homotopy left adjoint to the functor which sends a space to the constant $\infty$-sheaf on that space.
Moreover, $\sB$ is homotopy left adjoint to the inclusion of the $\infty$-category of constant $\infty$-sheaves on $\man$ into $\Sh(\man)$.
These two statements also hold if $\man$ is replaced by $\bD$, in which case we can also formulate the adjunction using $\infty$-presheaves instead of $\infty$-sheaves.
In particular, the inclusion $\Delta\to\bD$ is a homotopy initial functor (hence also an initial functor),
i.e., the functor $\tB$ computes the homotopy colimit over $\bD^\op$.
\end{proposition}

\begin{proof}
By a constant $\infty$-sheaf, we mean the homotopy sheafification of a constant presheaf.
Every constant presheaf on $\bD$ is an $\infty$-sheaf, and for a space $K$ the canonical functor from the constant $\infty$-sheaf to the mapping space $\infty$-sheaf $\textup{const}_K \to \map(-,K)$ is an objectwise weak equivalence.
On the other hand, a constant presheaf on $\man$ is in general not an $\infty$-sheaf.
However, since every open cover can be refined by a good open cover, the homotopy sheafification of a presheaf on $\man$ is determined by its restriction to $\bD$.
Therefore $\textup{const}_K \to \map(-,K)$ is also a weak equivalence in $\Sh(\man)$.

If $F$ is a representable presheaf, represented by a manifold $M$, then $\tB F \simeq M$.
From this it follows that
$$\map(F, \textup{const}_K) \simeq \map(\tB F, K)$$
for $F$ a representable, and extending by colimits, the same is true for any presheaf $F$, and then for any $\infty$-sheaf (the mapping space on the left is computed in the $\infty$-category of $\infty$-sheaves or equivalently, since $\textup{const}_K$ is an $\infty$-sheaf, in the $\infty$-category of presheaves).

As for the second statement, we have by \cref{mainthm} that $\sB F$ is a constant $\infty$-sheaf on $\man$, for any $\infty$-sheaf $F$.
Then, by Yoneda,
$\map(\sB F, \textup{const}_K) \simeq \map(\tB F, K)$, which by the first part implies the statement.

The argument for $\bD$ remains valid for $\infty$-presheaves because the constant $\infty$-presheaf is an $\infty$-sheaf in this case.
This also implies the claim about homotopy initiality of $\Delta\to\bD$.
\end{proof}

In other words, $\tB F$ is the shape of $\Sh(\man)/F$ as in Lurie~\cite[Chapter 7]{Lurie.HTT}, or equivalently, the fundamental $\infty$-groupoid of $F$ in the sense of Schreiber~\cite[Section 3.4]{Schreiber}, and $\sB F$ is the shape modality of $F$ in the sense of Schreiber~\cite[Definition 3.4.4]{Schreiber}.
For a different proof, see also Pavlov \cite[Proposition~12.10]{Pavlov.Diffeo}.
The statement about the homotopy adjunction continues to hold if $\infty$-sheaves on $\man$ are replaced by $\infty$-presheaves on $\man$,
and $\sB$ precomposed with the associated $\infty$-sheaf functor is weakly equivalent to~$\sB$,
i.e., the shape of the associated $\infty$-sheaf of an $\infty$-presheaf of~$F$ can be computed as the shape of~$F$ (Pavlov \cite[Proposition~13.9]{Pavlov.Diffeo}).

The following formal consequence of \cref{mainthm} proved to be useful in applications.
Sati–Schreiber \cite[Theorem~3.3.53]{SatiSchreiber} proposed the name “smooth Oka principle” for \cref{prop:internal},
in analogy to the Oka principle in complex geometry, and also gave an alternative derivation of \cref{prop:internal} from \cref{mainthm}.

\begin{proposition}\label{prop:internal}
Let $F$ be an $\infty$-sheaf on $\man$.
There is an evaluation map
$$\sB \Hom(M,F) \rightarrow \RHom(\sB M, \sB F)$$
which is a natural weak equivalence of $\infty$-sheaves for every manifold $M$.
Here $\Hom$ denotes the internal hom of $\infty$-sheaves, whereas $\RHom$ is the derived internal hom ($\Hom(M,-)$ is automatically derived).
\end{proposition}

\begin{proof}
The left side $\sB\Hom(M,F)$ is a concordance invariant $\infty$-sheaf because $\sB$ lands in concordance invariant presheaves by \cref{cor:enriched}
and $\infty$-sheaves by \cref{sheafthm}.
The right side $\Hom(\sB M, \sB F)$ is a concordance invariant $\infty$-sheaf because the left Bousfield localization that
produces concordance invariant $\infty$-sheaves is a cartesian localization (since open covers are closed under products with a fixed manifold)
and derived internal homs in cartesian left Bousfield localizations preserve local objects.
Therefore, the map under consideration is a map between concordance invariant $\infty$-sheaves,
so it is a weak equivalence if and only if its evaluation at the point is a weak equivalence.
Evaluating at the point gives the map
$$(\sB F)(M)\to \Rmap(\sB M,\sB F)\simeq \Rmap(\tB M,\tB F)=\Rmap(\Sing M,\tB F),$$
which is a weak equivalence by \cref{mainthm}.
\end{proof}


\medskip
Another application of this work is a construction of classifying spaces of field theories.
This has recently been done by Grady–Pavlov~\cite{Grady.Pavlov}, see in particular Theorem~8.2.9 there.
Stolz and Teichner have conjectured that concordance classes of particular classes of field theories determine cohomology theories~\cite{Stolz.Teichner}.
By Brown representability, this conjecture requires concordance classes of field theories to define a representable functor.
In brief, they define field theories as functors out of a category of bordisms equipped with a smooth map to a fixed manifold~$M$.
When the relevant bordism category is fully extended, field theories are an $\infty$-sheaf evaluated on~$M$.
The main result of this paper then shows that concordance classes of fully extended field theories are representable.
Furthermore, we identify a formula for the classifying space of field theories.

\subsection*{Acknowledgments}
We would like to thank Peter Teichner for asking the question that started this project in Spring 2011, and for his encouragement and patient feedback throughout the process.
We also thank Michael Weiss for generously sharing his insight into this problem (particularly in regard to \cref{lem:MW}).
This project has been going on for many years, and during that time we have benefited from discussions with many colleagues and hosting of various institutions.
In particular, we thank Ulrich Bunke, Jacob Lurie, Gustavo Granja, Aaron Mazel-Gee and Thomas Nikolaus for helpful discussions and comments, and the referee for helpful remarks that simplified an argument in a previous version of this paper.
We are grateful to the MSRI, the MPIM and the HIM in Bonn, and the universities of Münster, Göttingen, Regensburg, and Stanford for their hospitality.
D.P.'s gratitude extends to the $n$Lab for being a wonderful resource.
P.B.~was supported by FCT through grant SFRH/BPD/99841/2014.

\section{The concordance resolution is concordance-invariant}\label{sec:concinv}

\begin{notation}
Throughout, \emph{space} will mean simplicial set.
The category of such is denoted by $\Sp$.
A \emph{simplicial space} is a simplicial object in spaces, and the category of such is denoted by~$\sSp$.
Of course, this is the same as a bisimplicial set, though the terminology emphasizes that there is a preferred simplicial direction.
A simplicial set is often viewed as a simplicial \emph{discrete} space, by regarding a set as a discrete (or constant) simplicial set.
We denote the diagonal of a bisimplicial set~$X$ by~$|X|$,
this is our preferred model for homotopy colimits of simplicial spaces.
\end{notation}

\begin{definition}
\label{horns}
We write $\Delta^n$ for the representable simplicial \emph{set} and $\Delta[n]$ for the corresponding simplicial discrete space.
Similarly, we write $\partial \Delta^n$ and $\Lambda^n_k$ for the simplicial \emph{set} boundary and $k$th horn respectively,
and $\partial \Delta[n]$ and $\Lambda_k[n]$ for the corresponding simplicial spaces.
\end{definition}

\begin{definition}
We denote by $\man$
the (discrete) category of smooth manifolds (of any dimension) and smooth maps,
equipped with the Grothendieck topology of open covers.
\end{definition}

\begin{definition}
A presheaf $F: \man^{\op}\to \Sp$ is \emph{concordance invariant}
if for all manifolds $M$ the map $F(M)\to F(M\times \R)$ induced by the projection $M\times\R\to M$ is a weak equivalence.
\end{definition}

\begin{definition}
\label{defn:concres.man}
Set $$\A: \Delta\to\man,\qquad \A^n=\Bigl\{x\in\R^{n+1}\Bigm| \sum_i x_i=1\Bigr\}.$$
Given a presheaf $F: \man^{\op}\to \Sp$,
denote by $\sB F$ the presheaf
$$\sB F:\man^\op \to \Sp,\qquad \sB F(M) \coloneq \bigl|[k] \mapsto F(M \times \A^k)\bigr|,$$
where $|{-}|$ denotes the diagonal of a bisimplicial set.
\end{definition}

In this section we show that $\sB F$ is always concordance invariant.
Furthermore, if $\sB F$ is an $\infty$-\emph{sheaf}, then it is representable.
The arguments are largely formal,
so to make this structure more transparent we begin the discussion for an arbitrary category enriched over~$\Sp$ and later specialize to the category of manifolds.
These results are mostly a repackaging of Morel--Voevodsky~\cite{Morel.Voevodsky}, see also Herrmann--Strunk~\cite{Herrmann.Strunk}.

\begin{notation}
\label{defn:dcat}
Let $\sD_0$ be a discrete category with products and $\A^\bullet$ a cosimplicial object in $\sD_0$.
This data defines a category enriched in spaces, denoted by $\sD$, by declaring the set of $n$-simplices of the mapping space $\map_\sD(X,Y)$ to be $\hom_{\sD_0}(X \times \A^n, Y)$.
The example that will be of interest to us here is $\sD_0 = \man$.
\end{notation}

In this setting it makes sense to talk about concordance invariant presheaves on~$\sD_0$.

\begin{definition}\label{defn:conc}
Given a category $\sD_0$ as in \cref{defn:dcat},
a presheaf $F : \sD_0^\op \to \Sp$ is \emph{concordance invariant} if the map induced by the projection
$$F(X) \rightarrow F(X \times \A^1)$$
is a weak equivalence of spaces for all $X\in \sD_0$.
\end{definition}

A functor on $\sD_0$ that can be enriched, i.e., lifted to a functor on $\sD$,
necessarily sends smooth homotopies to simplicial homotopies and smooth homotopy equivalences to simplicial homotopy equivalences.
As such, it is automatically concordance invariant.
In \cref{prop:fibrant} below, we will prove the converse.

\begin{definition}\label{defn:concres}
Given a category $\sD_0$ as in \cref{defn:dcat},
the \emph{concordance resolution} of a functor $F : (\sD_0)^{\op} \to \Sp$
is the functor $$F(-\times {\A^{\bullet}}) : (\sD_0)^{\op} \to \sSp, \qquad X \mapsto F(X \times \A^\bullet).$$
We denote the homotopy colimit of $F( - \times \A^{\bullet})$ by
$$\sB F(X) \coloneq \hocolim_{[n] \in \Delta^\op} F(X \times \A^{n}) = \bigl| F(X\times \A^\bullet) \bigr|.$$
Here we use the diagonal of a bisimplicial set as a model for the homotopy colimit over $\Delta^\op$.
\end{definition}

For the case $\sD_0=\man$ this definition of $\sB F$ coincides with \cref{defn:concres.man}.

\begin{proposition}\label{prop:enriched}
Given a category~$\sD_0$ as in \cref{defn:dcat},
for any presheaf $F: \sD_0^\op \to \Sp$,
the functor $F(- \times \A^\bullet)$ lifts to an enriched functor $\sD^\op \to \sSp$.
\end{proposition}

\begin{proof}
We will define a simplicial map
$$\map(X,Y) \to \map(F(Y \times \A^{\bullet}), F(X \times \A^{\bullet})).$$
Let $g : X \times \A^n \to Y$ be a morphism in $\sD_0$.
Given a morphism $\alpha : [k] \to [n]$ in $\Delta$, consider the composition
$$F(Y \times \A^k) \xrightarrow{F(g \times \id_{\A^k})} F(X \times \A^n \times \A^k) \xrightarrow{F(\id_X \times \A^\alpha \times \id_{\A^k})} F(X\times\A^k\times\A^k) \xrightarrow{F(\id_X\times d)} F(X \times \A^k),$$
where $d: \A^k\to\A^k\times\A^k$ is the diagonal map.
This is functorial in $g$ and $\alpha$ and so defines a map between hom-sets
\begin{equation}\label{eq:enriched}
\hom(X \times \A^n, Y) \to \hom(F(Y \times \A^\bullet) \times \Delta[n], F(X \times \A^\bullet))
\end{equation}
for each $n \geq 0$.
Therefore $F(- \times {\A^\bullet})$ is enriched over spaces.
\end{proof}

\begin{corollary}\label{cor:enriched}
Given a category~$\sD_0$ as in \cref{defn:dcat},
for any presheaf $F: \sD_0^\op \to \Sp$,
the presheaf $\sB F$ is enriched over spaces and is concordance invariant.
\end{corollary}

\begin{proof}
To see that $\sB F$ is enriched, and hence concordance invariant, post-compose (\ref{eq:enriched}) with the homotopy colimit functor (alias geometric realization or diagonal) and use the fact that it commutes with products.
\end{proof}

\begin{remark}
The functor $\sB$ is homotopy left adjoint to the discretization functor
$$i^* : \PSh(\sD) \to \PSh(\sD_0),$$
given by the restriction along the inclusion $i:\sD_0\to\sD$.
This follows from the fact that $\sB$ is a left adjoint functor
whose value on a representable presheaf on $X\in\sD_0$
is the representable presheaf on $i(X)\in\sD$.
\end{remark}

The following proposition implies that the category of enriched presheaves on $\sD$ and the category of concordance-invariant presheaves on $\sD_0$ have equivalent homotopy theories.

\begin{proposition}\label{prop:fibrant}
Given a category~$\sD_0$ as in \cref{defn:dcat},
a presheaf $F: \sD_0^\op \to \Sp$ is concordance invariant if and only if the map $F(X)\to i^*\sB F(X)$ is a weak equivalence for all $X$.
\end{proposition}

\begin{proof}
If $F$ is concordance invariant then the simplicial object $F(X \times \A^{\bullet})$ is homotopically constant with value $F(X)$.
For the converse, consider the diagram
$$\begin{tikzpicture}[descr/.style={fill=white}]
\matrix(m)[matrix of math nodes, row sep=2.5em, column sep=2.5em,
text height=1.5ex, text depth=0.25ex] {
F(X)  & i^* \sB F(X) \\
F(X \times \A^1) & i^* \sB F (X \times \A^1) \\
};
\path[->] (m-1-1) edge node [auto] {} (m-1-2);
\path[->] (m-1-1) edge node [left] {} (m-2-1);
\path[->] (m-1-2) edge node [auto] {} (m-2-2) (m-2-1) edge node [auto] {} (m-2-2);
\end{tikzpicture}$$
The horizontal maps are weak equivalences by assumption.
The vertical map on the right is a weak equivalence since $\sB F$ is concordance invariant.
Thus, by the 2-out-of-3 property the left map is a weak equivalence.
\end{proof}

\begin{corollary}
\label{concord.presheaves}
The restriction functor
$$i^* : \PSh(\sD) \to \PSh(\sD_0)$$
is a right Quillen equivalence,
where $\PSh(\sD)$ is equipped with the projective model structure
and $\PSh(\sD_0)$ is equipped with the $\A^1$-invariant projective model structure,
i.e., the left Bousfield localization of the projective model structure with respect to the map $\A^1\to\A^0$.
\end{corollary}

\subsection{Concordance-invariant $\infty$-sheaves on manifolds are representable}

\begin{definition}\label{defn:sheaf}
A presheaf $F:\man^\op\to\Sp$ is an \emph{$\infty$-sheaf} if for every manifold $M$ and open cover $\{U_i \to M\}_{i \in I}$, the canonical map from $F(M)$ to the homotopy limit (over $\Delta$) of the cosimplicial space
$$\prod_{i_0 \in I} F(U_{i_0}) \leftleftarrows \prod_{i_0, i_1 \in I} F(U_{i_0} \cap U_{i_1}) \leftleftleftarrows \cdots$$
is a weak equivalence of spaces.
\end{definition}

Any $\infty$-sheaf $F$ satisfies $F(\varnothing) \simeq *$.
This is implied by the descent condition for the empty cover of the empty manifold.

\begin{remark}\label{rem:sheaves}
A set-valued sheaf is an $\infty$-sheaf of sets, and conversely.
Indeed, the (homotopy) limit of a cosimplicial discrete space is, by initiality, computed by the limit of its truncation to its $1$-coskeleton.
A \emph{stack} \emph{is} a groupoid-valued $\infty$-sheaf (Hollander~\cite{Hollander}).
Common alternative terminologies for $\infty$-sheaves include $\infty$-stacks and homotopy sheaves.
\end{remark}

The following proposition is due to Morel--Voevodsky~\cite{Morel.Voevodsky} and Dugger~\cite{Dugger}.

\begin{proposition}\label{prop:Dugger}
Given a presheaf $F: \man^\op\to\Sp$,
the presheaf $\sB F$ is an $\infty$-sheaf if and only if the evaluation map
$$\sB F(M) \to \map(\map(*,M), \tB F)$$
is a weak equivalence for every $M$,
where the evaluation map is the adjoint to the simplicial map
$$\map(*, M) \to \map(\sB F(M), \sB F(*)) \to \map(\sB F(M), \tB F)$$
gotten by the enrichment afforded by \cref{cor:enriched} and the map $\sB F(*)\to \tB F$ being the Kan complex replacement of $\sB F(*)$.
Here $\map(*, M)$ is weakly equivalent to $\Sing M$, the singular simplicial set of~$M$.
\end{proposition}

\begin{proof}
Take a good open cover $\{U_i \}_{i \in I}$ of $M$ and let $\sU_\bullet: \Delta^\op\to\man$ denote its Čech nerve.
There is a commutative square:
$$\begin{tikzpicture}[descr/.style={fill=white}]
\matrix(m)[matrix of math nodes, row sep=2.5em, column sep=2.5em,
text height=1.5ex, text depth=0.25ex] {
\sB F(M) & \Rmap(\Sing M, \sB F(*)) \\
\holim_{n\in\Delta} \sB F(\sU_n) & \holim_{n\in\Delta} \Rmap(\Sing \sU_n, \sB F(*)), \\
};
\path[->,font=\scriptsize] (m-1-1) edge node [auto] {} (m-1-2);
\path[->,font=\scriptsize] (m-1-1) edge node [left] {} (m-2-1);
\path[->, font=\scriptsize] (m-1-2) edge node [auto] {} (m-2-2) (m-2-1) edge node [auto] {} (m-2-2);
\end{tikzpicture}$$
where $\Sing$ denotes the singular simplicial set functor.
The right-hand vertical arrow is an equivalence since $\hocolim\nolimits_{n\in\Delta} \sU_n \simeq M$.
The lower horizontal arrow is an equivalence since $\sB F(V) \simeq \sB F(*)$ for $V$ contractible (by concordance invariance of $\sB F$).
The statement now follows by the $2$-out-of-$3$ property.
\end{proof}

\subsection{Homotopy groups of $\tB F$}

In this section, we explain how to compute the homotopy groups of $\tB F$. For a base-point $b$ in the $d$-dimensional sphere $S^d$ and $x\in F(*)$, let $\sB F(S^d)_x$ denote the homotopy fiber of $\sB F(b): \sB F(S^d)\to \sB F(*)$ over the image of~$x$ in~$\sB F(*)$.

\begin{proposition}
\label{homotopy.groups.bf}
Let $F:\man^\op\to\Sp$ be a presheaf satisfying the $\infty$-sheaf property with respect to finite covers, and let $x\in F(*)$.
The map
$$\pi_0 \sB F(S^d)_x \to \pi_d (\sB F(*),x)=\pi_d(\tB F, x)$$
is an isomorphism.

\end{proposition}

\begin{proof}
Under the assumption on $F$, the conclusion of \cref{mainthm} holds for any \emph{compact} manifold $M$.
This will be shown in \cref{cor:compact}.
Therefore, the top map in the commutative square
$$\begin{tikzpicture}[descr/.style={fill=white}]
\matrix(m)[matrix of math nodes, row sep=2.5em, column sep=2.5em,
text height=1.5ex, text depth=0.25ex] {
\sB F(S^d) & \Rmap(S^d, \sB F(*)) \\
\sB F(*) & \Rmap(*, \sB F(*)) \\
};
\path[->,font=\scriptsize] (m-1-1) edge node [auto] {} (m-1-2);
\path[->,font=\scriptsize] (m-1-1) edge node [left] {$\sB F(b)$} (m-2-1);
\path[->, font=\scriptsize] (m-1-2) edge node [auto] {$\Rmap(b, \sB F(*))$} (m-2-2) (m-2-1) edge node [auto] {} (m-2-2);
\end{tikzpicture}$$
is a weak equivalence.
The bottom map is a weak equivalence by construction.
Thus the induced map of vertical homotopy fibers over a point $x \in F(*)$ is a weak equivalence.
Taking $\pi_0$ of the map between homotopy fibers, we obtain the result.
\end{proof}

\begin{remark}
Elements in $\pi_0 \sB F(S^d)_x$ are concordance classes of sections of $F$ over $S^d$ which restrict to $x$ on $b \in S^d$.
We postpone the explanation to \cref{lem:pi0-sphere}.
\end{remark}

\begin{remark}
In the special case when $F$ is a concrete sheaf of sets, i.e., a diffeological space in the sense of Souriau~\cite{Souriau},
\cref{homotopy.groups.bf} resolves in the affirmative a conjecture of Christensen–Wu \cite[§1]{CW} on the isomorphism
of smooth homotopy groups of a diffeological space with the simplicial homotopy groups of its smooth singular complex.
Christensen–Wu \cite[Theorem~4.11]{CW} proved this conjecture for projectively fibrant diffeological spaces,
whereas \cref{homotopy.groups.bf} proves it for arbitrary simplicial presheaves satisfying the $\infty$-sheaf property for finite covers,
a much bigger class that includes all sheaves of sets, in particular, all diffeological spaces.
\end{remark}

\section{Weak Kan fibrations}\label{sec:weakKanfib}

As a warm-up to the ideas in this section, we will prove that the concordance relation~$\sim_{c}$ is an equivalence relation when $F$ is an $\infty$-sheaf.
This generalizes the standard fact that smooth homotopy is an equivalence relation, but the core of the argument is identical: gluing a pair of smooth maps along an open submanifold yields a smooth map.

\begin{lemma}\label{lem:concisequiv}
If $F: \man^\op\to\Sp$ is an $\infty$-sheaf, then $\sim_{c}$ is an equivalence relation on $F(M)_0$,
the set of 0-simplices in $F(M)$.
\end{lemma}

\begin{proof}
Reflexivity and symmetry are obvious.
To establish transitivity, suppose $\sigma_0$, $\sigma_1$ and $\sigma_2$ are such that $\sigma_0 \sim_{c} \sigma_1$ and $\sigma_1 \sim_{c} \sigma_2$.
Let $i_k$ denote the inclusion of $M \times\{k\}$ into $M \times \A^1$ and pick sections $\sigma_{01}$ and $\sigma_{12}$ over $M \times \A^1$ such that $i_0^*\sigma_{01} = \sigma_0$, $i_1^*\sigma_{01} = i_0^*\sigma_{12} = \sigma_1$ and $i_1^*\sigma_{12} = \sigma_2$.
Take a smooth map $r : \A^1 \to \A^1$ which fixes $0$ and $1$ and maps the complement of a small neighborhood of $1/2$ to $\{0,1\}$.
The sections $r^*\sigma_{01}$ and $r^*\sigma_{12}$ over $M \times \A^1$ are such that the restriction of $r^*\sigma_{01}$ to an \emph{open} neighborhood of $[1, \infty)$ agrees with the restriction of $r^* \sigma_{12}$ to an \emph{open} neighborhood of $(-\infty,0]$.
So, using the sheaf property and reparametrizing, we may glue these sections to obtain a section $\sigma_{012}$ over $M \times \A^1$ with $i_0^*\sigma_{012} = \sigma_0$ and $i_1^*\sigma_{012} = \sigma_2$, i.e., $\sigma_0 \sim_c \sigma_2$.
\end{proof}

The fact that $\sim_c$ is an equivalence relation is a shadow of an important property possessed by the concordance resolution: it is a \emph{0-weak Kan complex} (\cref{defn:wKancomplex}).
Informally, the \emph{weak} nature can be seen in the proof of transitivity at the point where sections $\sigma_{01}$ and $\sigma_{12}$ are replaced by $r^*\sigma_{01}$ and $r^*\sigma_{12}$.
This step is essential since sections cannot be glued along closed sets.
The failure of gluing along closed sets also means that concordance resolution does not satisfy the usual Kan condition as it does not have the right lifting property with respect to $\Lambda^2_1 \to \Delta^2$.
Similar features of smooth geometry allow us to show that certain restriction maps for the concordance resolution have analogous weak fibrancy properties.
The key definition formalizing this property is that of a \emph{weak Kan fibration}.

\subsection{Kan fibrations and weak Kan fibrations of simplicial spaces}

In this section, we define and investigate weak Kan fibrations of simplicial spaces (or sets).
These generalize Kan fibrations and are related to (and inspired by) Dold fibrations of topological spaces~\cite{Dold}.
We refer the reader to \cref{sec:appendix} for background on simplicial spaces.

The following definition is discussed by Lurie \cite{Lurie.Simplicial}, \cite[Definition A.5.2.1]{Lurie.SAG}.
Our definition is essentially the same, except that we formulate it for Reedy fibrant simplicial spaces,
to avoid mentioning Reedy fibrant replacements.

\begin{definition}\label{defn:Kanfib}
Let $f : X \to Y$ be a Reedy fibration between Reedy fibrant simplicial spaces.
We say that $f$ is a \emph{Kan fibration} if it has the right lifting property with respect to all horn inclusions (\cref{horns}).
That is, for every solid square
$$\begin{tikzpicture}[descr/.style={fill=white}, baseline=(current bounding box.base)]
\matrix(m)[matrix of math nodes, row sep=2.5em, column sep=2.5em,
text height=1.5ex, text depth=0.25ex] {
\Lambda_k[n] & X \\
\Delta[n] & Y \\
};
\path[->,font=\scriptsize] (m-1-1) edge node [auto] {} (m-1-2);
\path[->,font=\scriptsize] (m-2-1) edge node [auto] {} (m-2-2);
\path[->,font=\scriptsize] (m-1-1) edge node [left] {} (m-2-1);
\path[->,font=\scriptsize] (m-1-2) edge node [auto] {$f$} (m-2-2);
\path[->,dashed,font=\scriptsize] (m-2-1) edge node [auto] {} (m-1-2);
\end{tikzpicture}$$
there is a lift as pictured, where $n \geq 1$ and $0 \leq k \leq n$.
Similarly, we say that $f$ is a \emph{trivial Kan fibration} if it has the right lifting property with respect to $\partial \Delta[n] \to \Delta[n]$ for all $n \geq 0$.
\end{definition}

Unfortunately, \cref{defn:Kanfib} is not applicable in the situations of interest to us.
In particular, none of the maps in the crucial \cref{prop:CisKan,prop:EPisKan,prop:SnisKan}
satisfy \cref{defn:Kanfib}, so a result like Lurie~\cite[Theorem A.5.4.1]{Lurie.SAG} is not applicable.
An explicit counterexample is provided by \cref{Kan.counterexample}.
Therefore, we relax \cref{defn:Kanfib} to \cref{defn:wKan}.
First, we define an appropriate generalization of the right lifting property.

\begin{definition}\label{defn:wRLP}
Let $f : X \to Y$ be a Reedy fibration between Reedy fibrant simplicial spaces.
We say that $f$ has the \emph{weak right lifting property} (weak RLP)
with respect to a map $i: A \hookrightarrow B$
(and $i$ has the weak LLP with respect to~$f$)
if for every commutative square
$$\begin{tikzpicture}[descr/.style={fill=white}, baseline=(current bounding box.base)]
\matrix(m)[matrix of math nodes, row sep=2.5em, column sep=2.5em,
text height=1.5ex, text depth=0.25ex] {
A & X \\
B & Y \\
};
\path[->,font=\scriptsize] (m-1-1) edge node [auto] {$\beta$} (m-1-2);
\path[->,font=\scriptsize] (m-2-1) edge node [auto] {$\alpha$} (m-2-2);
\path[->,font=\scriptsize] (m-1-1) edge node [left] {$i$} (m-2-1);
\path[->,font=\scriptsize] (m-1-2) edge node [auto] {$f$} (m-2-2);
\path[->,dashed,font=\scriptsize] (m-2-1) edge node [auto] {$\widetilde{\alpha}$} (m-1-2);
\end{tikzpicture}$$
there is a lift $\widetilde{\alpha}$ as pictured, making the lower triangle commute strictly and the upper triangle commute \emph{up to a specified vertical homotopy}.
Such a homotopy consists of a map of simplicial spaces
$$H : A \times \Delta[1] \to X$$
subject to the requirement that $H_0 = \beta$, $H_1 = \widetilde{\alpha} \circ i$, and $f\circ H = \alpha\circ i \circ \pi$, where $\pi$ denotes the projection of $A \times \Delta[1]$ onto $A$.
\end{definition}

\begin{remark}\label{rem:mapcyl}
It will be useful to have some reformulations of \cref{defn:wRLP}.
Under the Reedy fibrancy hypothesis above, the map $f : X \to Y$ has the weak RLP with respect to $i : A \to B$ if and only if for every commutative square in the background of
$$\begin{tikzpicture}[
back line/.style={solid},
cross line/.style={preaction={draw=white, -,line width=6pt}}]
\matrix(m)[matrix of math nodes, row sep=1.0em, column sep=1.0em,
text height=1.5ex, text depth=0.25ex] {
A & & X \\
& M(i) & \\
B & & Y \\
& B & \\
};
\path[->,font=\scriptsize]
	(m-1-1) edge node [auto] {} (m-1-3)
	(m-1-1) edge node [left] {$i$} (m-3-1)
	(m-1-3) edge node [auto] {$f$} (m-3-3)
	(m-3-1) edge [back line] (m-3-3)
	(m-1-1) edge node [auto] {} (m-2-2)
	(m-2-2) edge [cross line] (m-4-2)
	(m-4-2) edge node [auto] {} (m-3-3)
	(m-1-1) edge node [auto] {} (m-2-2);
\path[->,font=\scriptsize] (m-3-1) edge node [left] {} (m-4-2);
\path[->,dashed,font=\scriptsize] (m-2-2) edge node [auto] {} (m-1-3);
\end{tikzpicture}$$
there exists a map from the \emph{mapping cylinder} $M(i) = B \sqcup_{i} (A \times \Delta[1])$ to $X$ making the diagram commute strictly.
Here the map in the foreground $M(i) \to B$ collapses $A \times \Delta[1]$ to $A$ -- we denote it by $\pi(i)$.

To put it differently, consider the category $\sSp^{[1]}  = \textup{Fun}(0 \to 1, \sSp)$ whose objects are maps of simplicial spaces and morphisms are commutative squares.
Then the requirement above is that the map induced by precomposition
$$\map\left(\downmap{M(i)}{B}, \downmap{X}{Y}\right) \to \map\left(\downmap{A}{B}, \downmap{X}{Y}\right)$$
(with the square on the left in the diagram above) is surjective on $0$-simplices.
Here $\map$ denotes the space of morphisms in $\sSp^{[1]}$.
\end{remark}

The Reedy model structure appears in the previous definitions as an artifact that guarantees homotopy invariance with respect to degreewise weak homotopy equivalences of simplicial spaces.
But it is possible (and worthwhile) to formulate a more homotopical definition of the weak RLP.

\begin{definition}\label{defn:wRLPgeneral}
A map $f : X \to Y$ between arbitrary simplicial spaces satisfies the \emph{weak right lifting property} (weak RLP) with respect to a map $i : A \to B$ if
$$\Rmap(\pi(i), f) \to \Rmap(i, f)$$
is surjective on $\pi_0$, where $\Rmap$ refers to the derived mapping space computed in the category $\sSp^{[1]}$ with objectwise weak equivalences.
\end{definition}

\begin{remark}\label{rem:wRLPinvariance}
We emphasize the homotopy invariance properties of this definition: a map $f$ has the weak RLP with respect to a map $i$ if and only if it has the weak RLP with respect to any map (degreewise) weakly equivalent to $i$.
Also, if a map $f$ has the weak RLP with respect to $i$ then so does any map (degreewise) weakly equivalent to $f$.
\end{remark}

\begin{proposition}
A Reedy fibration $f$ satisfies the weak RLP in the sense of \cref{defn:wRLP} if and only if it satisfies the weak RLP in the sense of \cref{defn:wRLPgeneral}.
\end{proposition}
\begin{proof}
Equip $\sSp$ with the Reedy (= injective) model structure and, relative to it, also equip $\sSp^{[1]}$ with the injective model structure.
In this model structure on $\sSp^{[1]}$, all objects are cofibrant.
Cofibrations are morphisms which are objectwise Reedy cofibrations of simplicial spaces, i.e., degreewise injections.
Fibrant objects are Reedy fibrations between Reedy fibrant simplicial spaces.

The morphism $i \to \pi(i)$, i.e., the commutative diagram
$$\begin{tikzpicture}[descr/.style={fill=white}, baseline=(current bounding box.base)]
\matrix(m)[matrix of math nodes, row sep=2.5em, column sep=2.5em,
text height=1.5ex, text depth=0.25ex] {
A & M(i) \\
B & B \\
};
\path[->,font=\scriptsize]
	(m-1-1) edge node [auto] {} (m-1-2)
	(m-2-1) edge node [auto] {} (m-2-2)
	(m-1-1) edge node [left] {} (m-2-1)
	(m-1-2) edge node [auto] {} (m-2-2);
\end{tikzpicture}$$
is a cofibration since the horizontal maps are degreewise injections.
Also, the map $f$ is a fibrant object in $\sSp^{[1]}$.
Therefore, the induced map
$$\map(\pi(i), f) \to \map(i, f)$$
is a fibration between Kan complexes, which is weakly equivalent to
$$\Rmap(\pi(i), f) \to \Rmap(i, f).$$
The result now follows, since a fibration of simplicial sets is surjective on $\pi_0$ if and only if it is surjective on $0$-simplices.
\end{proof}

\begin{definition}\label{defn:wKan}
A map $X \to Y$ between simplicial spaces is a \emph{weak Kan fibration} if it has the weak right lifting property with respect to the maps
$$h_i:\Sd^i (\Lambda_k[n])  \hookrightarrow \Sd^i (\Delta[n])$$
for each $i \geq 0$, $n \geq 1$ and $0 \leq k \leq n$.
Here the functor $\Sd:\sSp\to\sSp$ is as in \cref{defn:ex}.
\end{definition}

The definition also makes sense for maps $X \to Y$ of simplicial sets by regarding them as simplicial discrete spaces.

Note that every simplicial space is a weak Kan complex, in the sense that the map $X \to *$ is a weak Kan fibration.
This may seem odd at first, but it makes sense in light of \cref{sec:wkrfib}, as every space is tautologically quasi-fibrant.
A more interesting variation is
\begin{definition}
\label{defn:wKancomplex}
Given $\ell \geq -1$, a $\ell$-\emph{weak Kan complex} is a simplicial space~$X$
such that the terminal map $X \to *$ is a weak Kan fibration in which the vertical homotopies preserve the $\ell$-skeleton
of the horn inclusions~$h_i$.
\end{definition}

So a $(-1)$-weak Kan complex is simply a weak Kan complex, i.e., a simplicial space, and a $\infty$-weak Kan complex is a Kan complex (in the usual sense, as in \cref{defn:Kanfib}).
For a $0$-weak Kan complex~$X$, the image of $\pi_0X_1 \to \pi_0 X_0 \times \pi_0 X_0$ is an equivalence relation (cf.~\cref{lem:concisequiv}).
We will not make use of the notion of $\ell$-weak Kan complex for $\ell > 0$, and for $\ell = 0$ we will use it in \cref{cor:infinite-prods}.

\begin{example}
A Kan fibration (\cref{defn:Kanfib}) is a weak Kan fibration.
Of course, the usual definition of a Kan fibration does not mention subdivisions;
this is because a map satisfying the strict RLP against all horn inclusions $\Lambda_k[n] \to \Delta[n]$
automatically satisfies the same property against all subdivisions of those,
since these subdivided horn inclusions can be presented using cobase changes and the strict LLP is stable under cobase change,
so the strict LLP for horn inclusions implies the strict LLP for subdivided horn inclusions.
The same is \emph{not} true for the weak RLP, so we need to take subdivisions seriously.

Likewise, if $X$ is a $\ell$-weak Kan complex and $K \hookrightarrow L$ is an inclusion of simplicial sets, it does not follow automatically that $X^L \to X^K$ is a weak Kan fibration.
\end{example}

\begin{lemma}\label{lem:ex}
Kan's $\Ex$ functor (see \cref{defn:ex}) preserves weak Kan fibrations and so does $\Ex^\infty$.
\end{lemma}
\begin{proof}
The functor $\Ex$ is right adjoint to the subdivision functor, so we are investigating a square of the form
\begin{equation}\label{eq:ex}
\begin{tikzpicture}[descr/.style={fill=white}, baseline=(current bounding box.base)]
\matrix(m)[matrix of math nodes, row sep=2.5em, column sep=2.5em,
text height=1.5ex, text depth=0.25ex] {
\Sd^{i+1} \Lambda_k[n] & X \\
\Sd^{i+1} \Delta[n] & Y \\
};
\path[->,font=\scriptsize] (m-1-1) edge node [auto] {} (m-1-2);
\path[->,font=\scriptsize] (m-2-1) edge node [auto] {} (m-2-2);
\path[->,font=\scriptsize] (m-1-1) edge node [left] {} (m-2-1);
\path[->,font=\scriptsize] (m-1-2) edge node [auto] {} (m-2-2);
\path[->,dashed,font=\scriptsize] (m-2-1) edge node [auto] {} (m-1-2);
\end{tikzpicture}
\end{equation}
Since $f$ is a weak Kan fibration, there is a lift as shown, and a homotopy $H : \Sd^{i+1} \Lambda_k[n] \times \Delta[1] \to X$.
Since there is always a map from the subdivision of the product to the product of the subdivisions, we can precompose $H$ with
$$\Sd(\Sd^i \Lambda_k[n] \times \Delta[1]) \to \Sd^{i+1} \Lambda_k[n] \times \Sd \Delta[1] \to \Sd^{i+1} \Lambda_k[n] \times \Delta[1].$$
This gives the required homotopy for the upper triangle which is vertical over $Y$, proving the result.
The case of $\Ex^\infty$ follows automatically since a map from a finite-dimensional simplicial space to $\Ex^\infty X$ factors through some finite $\Ex^i X$.
\end{proof}

\begin{remark}
Weak Kan fibrations are stable under various operations.
They are stable under \emph{homotopy} base change (with respect to degreewise weak equivalences of simplicial spaces).
In other words, weak Kan fibrations that are moreover Reedy fibrations between Reedy fibrant objects are stable under pullback.

Weak Kan fibrations are also stable under fiberwise homotopy retracts (that is, if $g : W \to Y$ is a homotopy retract over $Y$ of a weak Kan fibration $f : X \to Y$ then $g$ is a weak Kan fibration).
In particular, weak Kan fibrations are stable under fiberwise homotopy equivalences.
Moreover, if we allow subdivisions of the vertical homotopies in the definition of the weak lifting property,
i.e., if we replace $\Delta[1]$ by $\Sd^i \Delta[1]$ for $i \geq 0$ in \cref{defn:wRLP}, then a composition of two weak Kan fibrations is also a weak Kan fibration.
(Although, the resulting notion would presumably be weaker than the one we are using.)
Since these properties will not be used in what follows, and the proofs are not particularly difficult, we omit further explanations.
\end{remark}

\subsection{Weak Kan fibrations are realization fibrations}\label{sec:wkrfib}

\begin{definition}[Rezk~\cite{Rezk.Commute}]\label{defn:realfib}
A map $f : X \to Y$ of simplicial spaces is a \emph{realization fibration} if for every $Z \to Y$ the induced map
\begin{equation}\label{eq:realfib}
| X \times^h_{Y} Z | \to |X| \times^h_{|Y|} |Z|
\end{equation}
is a weak equivalence of spaces.
The vertical bars refer to the diagonal simplicial set, which models the homotopy colimit over $\Delta^\op$.
\end{definition}

\begin{remark}
Realization fibrations are related to quasi-fibrations in the sense of Dold--Thom~\cite{Dold.Thom}.
For example, if a map $f : X \to Y$ between simplicial \emph{sets} is a realization fibration, then the map (\ref{eq:realfib}) with $Z$ a point is identified with the inclusion of the fiber of $f$ into the homotopy fiber.
That is, $f$ is a quasi-fibration.
On the other hand, not all quasi-fibrations are realization fibrations: realization fibrations are stable under homotopy pullback, whereas quasi-fibrations need not be.
\end{remark}

We now turn to the main technical result of the section.
\begin{theorem}\label{thm:wkanmain}
A weak Kan fibration is a realization fibration.
\end{theorem}

As we already emphasized, the subdivisions of simplices and horns in the definition of a weak Kan fibration are important for a number of reasons.
\cref{lem:ex} is one such reason, that will be exploited later on.
Below is another.
\begin{example}
This is an example of a map which has the weak right lifting property against the (nonsubdivided) map $(h_0)$ but is not a realization fibration.
Suppose $X$ is the union of three nondegenerate $1$-simplices as in the picture:
$$\begin{tikzpicture}[descr/.style={fill=white}, baseline=(current bounding box.base)]
\matrix(m)[matrix of math nodes, row sep=1em, column sep=1em,
text height=1.0ex, text depth=0ex] {
\bullet & \bullet \\
\bullet & \bullet \\
};
\path[->,font=\scriptsize] (m-1-1) edge node [auto] {$$} (m-1-2);
\path[->,font=\scriptsize] (m-2-1) edge node [auto] {$$} (m-2-2);
\path[->,font=\scriptsize] (m-1-1) edge node [auto] {} (m-2-1);
\end{tikzpicture}$$
and $Y = \Delta[1]$.
Let $f : X \to Y$ be the projection in the vertical direction.
This is not a Kan fibration: there are $2$-horns $\Lambda^k[2]$ in $X$ that cannot be filled.
On the other hand, $f$ has the \emph{weak} right lifting property with respect to the map $(h_0)$.
But $f$ is not a quasi-fibration, since at one point the fiber is disconnected while the homotopy fiber is not.
Hence $f$ cannot be a realization fibration.
This is not a contradiction: $f$ is not a weak Kan fibration as it does not have the weak right lifting property for $(h_2)$, the second subdivision of the horn inclusion.
\end{example}

\begin{proposition}\label{prop:weakGJ}
Suppose $f : X \to Y$ is a Reedy fibration between Reedy fibrant simplicial spaces.
Then $f$ is a weak Kan fibration of simplicial spaces if and only if $|f| : |X| \to |Y|$ is a weak Kan fibration of simplicial sets.
The vertical bars denote the diagonal functor.
\end{proposition}

\begin{proof}
Let us denote by $\kappa$ the map $\sd^i \Lambda^n_k \to \sd^i \Delta^n$ for $i \geq 0$.
Write $\delta$ for the functor which sends a simplicial set $K$ to the corresponding simplicial discrete space $[n] \mapsto K_n$.
Recall the commutative square of spaces $\kappa \to \pi(\kappa)$ where $\pi(\kappa)$ is the projection of mapping cylinder $M(\kappa)$ onto $\sd^i \Delta^n$.
By \cref{defn:wKan}, the map $f$ is a weak Kan fibration if and only if
$$\Rmap(\delta(\pi(\kappa)), f ) \to \Rmap(\delta(\kappa), f )$$
is surjective on $\pi_0$.

The realization (i.e., diagonal) functor has a left adjoint $d_! : \Sp \to \sSp$, and there is a natural transformation $d_! \to \delta$ which is a weak equivalence (\cref{lem:d!}).
Therefore, $f$ is a weak Kan fibration if and only if
$$\Rmap(d_! \pi(\kappa), f ) \to \Rmap(d_! \kappa, f ).$$
is surjective on $\pi_0$.
Regarding $\sSp^{[1]}$ with the injective model structure, the map $d_! \kappa \to d_! \pi(\kappa)$ is a cofibration between cofibrant objects (since $d_!$ sends monomorphisms to monomorphisms), and the target~$f$ is fibrant by hypothesis.
As such, the above holds if and only if the Kan fibration between Kan simplicial sets
$$\map(d_! \pi(\kappa), f ) \to \map(d_!(\kappa), f )$$
is surjective.
By adjunction, this holds if and only if $|f|$ has the weak RLP with respect to $\kappa$.
(Note that $|f|$, being a map between simplicial discrete spaces, is automatically a Reedy fibration between Reedy fibrant objects, and so a fibrant object in $\sSp^{[1]}$.)
\end{proof}

\begin{remark}\label{rem:strongKan}
The same proof, with $\kappa$ of the form $\Lambda^n_k \to \Delta^n$ and with $\pi(\kappa)$ replaced by the identity $\Delta^n \to \Delta^n$,
shows that a Reedy fibration between Reedy fibrant simplicial spaces $f$ is a (trivial) Kan fibration if and only if $|f|$ is a (trivial) Kan fibration.
\end{remark}

In order to prove that weak Kan fibrations are realizations fibrations, we will use the following criterion due to Rezk.

\begin{theorem}[Rezk~\cite{Rezk.Commute}]\label{thm:Rezk}
A map $f : X \to Y$ of simplicial spaces is a realization fibration if and only if for all maps $\Delta[m] \to Y$ and $\Delta[0] \to \Delta[m]$, the induced map on pullbacks
$$X \times^h_Y \Delta[0] \to X \times^h_Y \Delta[m]$$
is a weak equivalence after realization.
\end{theorem}

\begin{proof}[Proof of \cref{thm:wkanmain}]
Let $f : X \to Y$ be a weak Kan fibration.
We will verify that $f$ satisfies the condition in \cref{thm:Rezk}.
We may assume, without loss of generality, that $f$ is a Reedy fibration between Reedy fibrant simplicial spaces.
Then the homotopy pullbacks above become pullbacks.
Using the fact that realization and $\Ex^\infty$ commute with finite limits, our task is then to show that
$$\Ex^\infty |X| \times_{\Ex^\infty |Y|} \Ex^\infty \Delta^0 \to \Ex^\infty |X| \times_{\Ex^\infty |Y|} \Ex^\infty \Delta^m$$
is a weak equivalence of Kan complexes.
Since $f$ is a weak Kan fibration, the same can be said of $|f|$ by \cref{prop:weakGJ} and of $\Ex^\infty |f|$ by \cref{lem:ex}.
To simplify the notation, let us write $g : U \to V$ for $\Ex^\infty |f| : \Ex^\infty |X| \to \Ex^\infty |Y|$.

In view of \cref{prop:duggerisaksen} (and \cref{cor:duggerisaksen} and \cref{expl:duggerisaksen}), we will show that for every solid diagram
$$\begin{tikzpicture}[
back line/.style={solid},
cross line/.style={preaction={draw=white, -,line width=6pt}}]
\matrix(m)[matrix of math nodes, row sep=1.0em, column sep=1.0em,
text height=1.5ex, text depth=0.25ex] {
\partial \Delta^n & & U \times_{V} \Delta^0 \\
& \partial \Delta^n \times \Delta^1 \sqcup_{\partial \Delta^n \times \{1\}} \Delta^n & \\
\Lambda^{n+1} &  & U \times_{V} \Ex^\infty \Delta^m \\
& \Lambda^{n+1} \times \Delta^1 \sqcup_{\Lambda^{n+1}  \times \{1\}} \Delta^{n+1} & \\
};
\path[->,font=\scriptsize]
	(m-1-1) edge node [auto] {} (m-1-3)
	(m-1-1) edge node [left] {$i$} (m-3-1)
	(m-1-3) edge node [auto] {} (m-3-3)
	(m-3-1) edge [back line] (m-3-3)
	(m-1-1) edge node [auto] {} (m-2-2)
	(m-2-2) edge [cross line] (m-4-2)
	(m-1-1) edge node [auto] {} (m-2-2)
	(m-3-1) edge node [auto] {} (m-4-2);
\path[->,dashed,font=\scriptsize]
	(m-4-2) edge node [auto] {} (m-3-3)
	(m-2-2) edge node [auto] {} (m-1-3);
\end{tikzpicture}$$
there are dashed maps as pictured.
Let us write $A \to B$ for the middle vertical arrow.
Consider the map
\begin{equation}\label{eq:bla2}
\Lambda^{n+1} \sqcup_{\partial \Delta^{n} \times \{0\}} A \to \Ex^\infty \Delta^m
\end{equation}
determined by the lower horizontal map in the diagram above and by the map $A \to \Delta^0 \to \Ex^\infty \Delta^m$.
Since the terminal map $\Ex^\infty \Delta^m \to \Delta^0$ is a trivial Kan fibration, the map (\ref{eq:bla2}) extends along the inclusion
$$\Lambda^{n+1} \sqcup_{\partial \Delta^{n} \times \{0\}} A \hookrightarrow B.$$
Next, we want to define a map $B \to U$ which is compatible with the composition
$$B \to \Ex^\infty \Delta^m \to V$$ of the map we have just constructed.
This will give us the lower dashed map in the diagram above.
Consider the diagram
$$\begin{tikzpicture}[descr/.style={fill=white}, baseline=(current bounding box.base)]
\matrix(m)[matrix of math nodes, row sep=2.5em, column sep=2.5em,
text height=1.5ex, text depth=0.25ex] {
\Lambda^{n+1} & U \\
\Delta^{n+1} & V \\
};
\path[->,font=\scriptsize]
	(m-1-1) edge node [auto] {} (m-1-2)
	(m-2-1) edge node [auto] {} (m-2-2)
	(m-1-1) edge node [left] {} (m-2-1)
	(m-1-2) edge node [auto] {$g$} (m-2-2);
\end{tikzpicture}$$
where the lower map is the composition $\Delta^{n+1} \to B \to \Ex^\infty \Delta^m \to V$.
Since $g$ is a weak Kan fibration, we obtain a lift and vertical homotopy, i.e., the required map $B \to U$.
Therefore, we have defined a map $B \to U \times_V \Ex^\infty \Delta^m$, whose restriction to $A$ factors through $U \times_V \Delta^0$.
\end{proof}

\section{Weak Kan fibrancy of the concordance resolution}\label{sec:weakKanconc}

Let $F:\man^\op\to\Sp$ be an $\infty$-sheaf.
In this section and the next, we apply the general theory developed in the previous sections to prove \cref{sheafthm}.
The goal of this section is to prove the following three propositions.

\begin{proposition}\label{prop:CisKan}
For any $\infty$-sheaf $F:\man^\op\to\Sp$ the simplicial space $[n] \to F(\A^n)$ is a $0$-weak Kan complex.
\end{proposition}

\begin{proposition}\label{prop:EPisKan}
For any $\infty$-sheaf $F:\man^\op\to\Sp$ the map of simplicial spaces
$$F(\A^\bullet \times \A^1 ) \to F(\A^\bullet \times \partial \A^1)$$
is a weak Kan fibration.
Here $\partial \A^1=\A^0\sqcup \A^0$ is the disjoint union of two points.
\end{proposition}

\begin{proposition}\label{prop:SnisKan}
Let $*\hookrightarrow S^n$ be an inclusion of a basepoint into the smooth $n$-dimensional sphere.
For any $\infty$-sheaf $F:\man^\op\to\Sp$,
the induced map of simplicial spaces
$$F(\A^\bullet \times S^n)\to F(\A^\bullet)$$
is a weak Kan fibration.
\end{proposition}

We use these to prove \cref{sheafthm} in \cref{sec:sheafprop}, and the reader may wish to jump directly to that section to see these propositions in action.
In fact, in that section and in the remainder of the paper, we will only need \cref{prop:EPisKan,prop:SnisKan}.
We include \cref{prop:CisKan} because its proof anticipates much of the proof of the other two propositions.

The proofs of these propositions are based on the following simple observation.
If one were to try to prove that $F(\A^\bullet)$ is a Kan complex in the sense of \cref{defn:Kanfib},
an obvious approach would be to construct a deformation retraction from the simplex~$\A^n$ to its horn.
This is possible in the topological and PL settings, but not possible smoothly,
and a rigorous proof of this is provided in \cref{Kan.counterexample}.
However, the basic idea can be salvaged if one only asks for~$F(\A^\bullet)$ to be $0$-\emph{weak} Kan, which roughly translates into asking for a retraction up to a suitable homotopy.
This parallels the proof that concordance is an equivalence relation (\cref{lem:concisequiv}): we modify smooth maps between manifolds (via homotopies) to achieve certain constancy properties.
A relative version of this line of reasoning applies to the maps in \cref{prop:EPisKan,prop:SnisKan}.

\subsection{The sheaf associated to a simplicial set}

\begin{definition}
\label{adjunction.pre}
Denote by $$\|{-}\|_\pre: \sSp \to \PSh(\man,\Sp)$$
the simplicial left adjoint functor that sends $\Delta[n]$ to the representable presheaf of $\A^n$, $\hom(-,\A^n)$.
The corresponding simplicial right adjoint functor is $$\PSh(\man,\Sp)\to\sSp, \qquad F \mapsto (n\mapsto F(\A^n)).$$
\end{definition}

For a simplicial space $K$, the presheaf $\|K\|_\pre: \man^\op \to \Sp$ is given by the coend
$$K_n \otimes_{[n] \in \Delta} \hom(-,\A^n).$$
If $K$ is a simplicial set (i.e., a simplicial discrete space), then $\|K\|_\pre$ is a presheaf of sets.

\begin{remark}
The adjunction of \cref{adjunction.pre} is Quillen if both categories are equipped with the projective model structure or with the injective model structure.
Therefore, there is a weak equivalence
$$\Rmap(\|K\|_{\pre}, F) \xrightarrow{} \Rmap(K,F(\A^\bullet))$$
natural in the simplicial space $K$ and the presheaf~$F$.
\end{remark}

\begin{remark}
Note that $\| K \|_{\pre}$ is usually not a sheaf of sets.
For a simple example illustrating this, take $\| \Lambda^2_1 \|$ which is the pushout $$\hom(-,\A^1) \sqcup_{\hom(-,\A^0)} \hom(-,\A^1)$$ and pick an open cover of $\R^1$ by two open sets and compatible sections over each that do not lift to a section over~$\R^1$.
\end{remark}

\begin{definition}
Given a simplicial set~$K$, denote by $\|K\|$ the associated sheaf of sets of $\|K\|_\pre$.
\end{definition}

Recall that for presheaves of \emph{sets}, the notions of a sheaf and $\infty$-sheaf coincide (\cref{rem:sheaves}).
For such presheaves, sheafification and $\infty$-sheafification also agree.
A reference is Dugger–Hollander–Isaksen \cite[Proposition A.2]{Dugger.Hollander.Isaksen}.
(One can also prove this directly by comparing the usual sheafification formula, given by the so-called plus construction,
to its $\infty$-counterpart applied to a presheaf of sets.)
If we tried to define $\|K\|$ for simplicial \emph{spaces}~$K$, we would have had to use $\infty$-sheafification from the start.

\begin{remark}
Since sheafification is left adjoint to the inclusion of sheaves into presheaves, for $F$ an $\infty$-sheaf, the weak equivalence above lifts to a weak equivalence
$$\Rmap(\|K\|, F) \xrightarrow{\simeq} \Rmap(K,F(\A^\bullet))$$
natural in the simplicial space $K$ and the $\infty$-sheaf $F$.
\end{remark}

\begin{remark}
To be more concrete, suppose $K$ is the simplicial set associated to a simplicial complex with vertex set $V$ and take $X$ to the union of affine subspaces spanned by the simplices of $K$.
Then $||K||(M)$ is the set of smooth maps $M \to \R^V$ that land in $X$.
\end{remark}

\begin{example}
\label{Kan.counterexample}
Consider the sheaf of sets $F=\|\Lambda_1[2]\|$.
The simplicial object $n\mapsto F(\A^n)$ is not a Kan complex in the sense of \cref{defn:Kanfib}.
Indeed, consider the horn $\Lambda_1[2]\to (n\mapsto F(\A^n))$ given by the adjoint map of the identity map on~$F$.
This horn does not admit a filling by~$\Delta[2]$.
Indeed, such a filling has to be a section $s\in F(\A^2)$ that restricts to the identity map on $\|\Lambda_1[2]\|$.
Since $\|\Lambda_1[2]\|$ is the sheafification of the sheaf of sets $\|\Lambda_1[2]\|_\pre$,
in some neighborhood of the vertex $1\in \A^2$, the section~$s$ must factor through one of the two 1-dimensional faces of $\Lambda_1[2]$.
However, possessing such a factorization means that $s$ cannot restrict to the identity map on any neighborhood of the vertex~1 in $\|\Lambda_1[2]\|$.
Thus, $n\mapsto F(\A^n)$ is not a Kan complex in the sense of \cref{defn:Kanfib}.
\end{example}

\begin{proposition}\label{prop:inj}
Let $F$ be an injectively fibrant object of $\PSh(\man,\Sp)$.
Then the maps in \cref{prop:EPisKan,prop:SnisKan} are Reedy fibrations between Reedy fibrant simplicial spaces.
\end{proposition}

\begin{proof}
We show that the map in \cref{prop:EPisKan} is a Reedy fibration.
The argument for the one in \cref{prop:SnisKan} is similar.
Let $A \to B$ be a trivial Reedy cofibration of simplicial spaces, i.e., a map of simplicial spaces which is a degreewise monomorphism and a degreewise weak equivalence.
Using the adjunction of \cref{adjunction.pre},
denote by~$Q$ the pushout of presheaves
\begin{equation}\label{eq:Q}
\begin{tikzpicture}[descr/.style={fill=white}, baseline=(current bounding box.base)]
\matrix(m)[matrix of math nodes, row sep=2.5em, column sep=2.5em,
text height=1.5ex, text depth=0.5ex] {
\| A \|_\pre \times \partial \A^1 & \| B \|_\pre \times \partial \A^1 \\
\| A \|_\pre  \times \A^1  & Q \\
};
\path[->,font=\scriptsize]
	(m-1-1) edge node [auto] {$$} (m-1-2)
	(m-1-1) edge node [left] {$$} (m-2-1)
	(m-1-2) edge node [auto] {$$} (m-2-2)
	(m-2-1) edge node [below] {$$} (m-2-2);
\end{tikzpicture}
\end{equation}
By adjunction, $F(\A^\bullet \times \A^1 ) \to F(\A^\bullet \times \partial \A^1)$ has the right lifting property with respect to $A \to B$ if and only if every solid diagram of presheaves
\begin{equation*}
\begin{tikzpicture}[descr/.style={fill=white}, baseline=(current bounding box.base)]
\matrix(m)[matrix of math nodes, row sep=2.5em, column sep=2.5em,
text height=1.5ex, text depth=0.5ex] {
Q & F \\
\| B \|_\pre \times \A^1 & \\
};
\path[->,font=\scriptsize] (m-1-1) edge node [auto] {$$} (m-1-2);
\path[->,font=\scriptsize] (m-1-1) edge node [left] {$$} (m-2-1);
\path[->,dashed,font=\scriptsize] (m-2-1) edge node [below] {$$} (m-1-2);
\end{tikzpicture}
\end{equation*}
has a lift as pictured.
The existence of this lift is part of the so-called pushout-product axiom, which holds for the category of presheaves on $\man$.
But we'll provide a short argument here.
The top horizontal map in square (\ref{eq:Q}) is a trivial cofibration since $\|{-}\|_\pre$ is a left Quillen functor.
Since trivial cofibrations are stable under cobase change, the lower horizontal map is also a trivial cofibration.
But the composition
$$\| A \|_\pre  \times \A^1 \to Q \to  \| B \|_\pre \times \A^1$$
is also a trivial cofibration by hypothesis, and so by two-out-of-three the right-hand map is a weak equivalence.
The right-hand map is also an injective cofibration, as can be checked directly.
Since $F$ is injectively fibrant, we conclude that the dashed map exists.
\end{proof}

\begin{lemma}\label{lem:pi0-sphere}
Let $F:\man^\op\to\Sp$ be a presheaf satisfying the $\infty$-sheaf property with respect to finite covers and let $x \in F(*)$.
Let $F(S^d)_x$ denote the homotopy fiber, over $x$, of the map $F(S^d) \to F(*)$ induced by a choice of basepoint in the $d$-sphere $S^d$.
There is a canonical isomorphism
$$\left({\pi_0 F(S^d)_x } \right) / {\sim} \to \pi_0 \sB F(S^d)_x$$
where $\sim$ is the equivalence relation of concordance and $\sB F(S^d)_x$ is the homotopy fiber of $\sB F(S^d) \to \sB F(*)$ as in \cref{homotopy.groups.bf}.
Therefore, the quotient of $\pi_0 F(S^d)_x$ by the equivalence relation of concordance is isomorphic to $\pi_d (\sB F(*), x)$.
\end{lemma}

\begin{proof}
By \cref{prop:SnisKan}, the map induced by a choice of basepoint in $S^d$
$$F(\A^\bullet \times S^d)\to F(\A^\bullet)$$
is a weak Kan fibration.
Applying the homotopy colimit functor to this map and then taking its homotopy fiber over $x \in F(\A^0)$ yields the space $\sB F(S^d)_x$, by definition.

Let $Z_\bullet \coloneq \textup{hofiber}_{x}  (F(\A^\bullet \times S^d) \to F(\A^\bullet))$.
The canonical map
$$|Z_\bullet | \to \sB F(S^d)_x$$
is a weak equivalence since weak Kan fibrations are realization fibrations (\cref{thm:wkanmain}).
Moreover, weak Kan fibrations are stable under base change, so $Z_\bullet$ is a $0$-weak Kan complex and as such
the relation on $\pi_0 Z_0$ determined by the faces $\pi_0 Z_1 \to \pi_0 Z_0 \times \pi_0 Z_0$ is an equivalence relation (cf.~\cref{lem:concisequiv}).
This equivalence relation is the concordance relation on $\pi_0 Z_0 = \pi_0 F(S^d)_x$.

The second statement follows from \cref{homotopy.groups.bf}.
\end{proof}

\subsection{Closed simplices}\label{sec:closedsimplices}

\begin{definition}
\label{bdelta}
Denote by $\bDelta^n:\man^\op\to\Sets$ the subsheaf of the representable sheaf of sets $\hom(-,\A^n)$
consisting of sections $X \to \A^n$ whose image is contained in the closed simplex $\R^{n+1}_{\ge0}\cap \A^n \subset \A^n$.
\end{definition}


\begin{definition}
Given a simplicial set~$K$,
denote by $\overline{\| K \|}:\man^\op\to\Sets$ the subsheaf of sets of ${\| K \|}$ consisting of smooth maps $U \to M$
that factor locally through some \emph{closed} simplex of $K$ in $M$.
In formulas,
$$\overline{\| K \|} = L(\overline{\| K \|}_\pre) = L(K_n \otimes_{[n] \in \Delta} \bDelta^n)$$
where $\bDelta^n$ is as in \cref{bdelta},
$L(W)$ denotes the associated sheaf of a presheaf of sets~$W$,
and $\overline{\| K \|}_\pre$ denotes the left adjoint
functor of $F \mapsto \map(\bDelta^\bullet,F)$ applied to the simplicial set~$K$.
\end{definition}

Concretely, if $K$ is a triangulation of a smooth manifold $N$, a section of $\| K \|$ over a manifold $M$ is a smooth map $M \to N$ which factors \emph{locally} through some simplex of the triangulation.


\begin{lemma}\label{lem:resheq}
Let $K$ be a simplicial set for which faces of nondegenerate simplices are nondegenerate.
Then the sheaf $\| K \|$ weak deformation retracts onto $\overline{\| K\|}$.
That is, there is a map $h : \| K \| \times \A^1 \to \| K\|$ whose restriction to $\|K \| \times \{0\}$ is the identity, whose restriction to $\| K \| \times \{1\}$ factors through $\overline{\| K\|}$, and $h( \overline{\|K\|} \times \A^1) \subset \overline{\|K\|}$.
Moreover, for any subcomplex $L$ of $K$, $h$ restricts to a weak deformation retraction of $\| L \|$ onto $\overline{\| L\|}$.
\end{lemma}

\begin{proof}
Consider a homotopy
$$\lambda^n : \A^n\times [0,1] \to \A^n$$
that in barycentric coordinates is constructed as follows.
Fix $c : \A^1 \to [0,\infty)$ a smooth function with $c \equiv 0$ on $(-\infty,0]$ and strictly increasing on $[0,\infty)$.
Then we set
$$\lambda^n((x_0, \ldots, x_n), t) = ( y_0/C_t, \ldots, y_n/C_t)$$
where $y_i = t c(x_i) + (1-t)x_i$ and $C_t = \sum_{i = 0}^n y_i$. 
Extend $\lambda^n$ to $\A^n \times \A^1$ by precomposing $\lambda^n$ with $\id \times f : \A^n \times \A^1 \to \A^n \times [0,1]$ where $f$ is a smooth function which takes value $0$ on a neighborhood of $(-\infty,0]$ and $1$ on a neighborhood of $[1, \infty)$.
This gives a weak deformation retraction $h : \A^n \times \A^1 \to \A^n$ of $\A^n$ onto the subsheaf $\bDelta^n$ for each fixed $n$.

The map $\lambda$ is functorial with respect to injections $[m] \to [n]$. And we may replace the category~$\Delta$ in the coend defining $\| K \|_\pre$ with the subcategory $\Delta_\inj$ of injective maps. To see this, we can express the coend over $\Delta$, respectively $\Delta_\inj$, as a colimit over the category $\textup{simp}(K)$ of simplices of $K$, respectively, the subcategory $\textup{ndsimp}(K)$ of nondegenerate simplices. The inclusion $\textup{ndsimp}(K) \to \textup{simp}(K)$ is terminal (alias final) under the assumption that each face of a nondegenerate simplex of $K$ is nondegenerate. Hence, $\lambda$ defines a weak deformation retraction $h: \|K\|_\pre \times \A^1 \to \|K\|_\pre$ of $\| K \|_\pre$ in $\overline{\| K \|}_\pre$.


Now, composing $h$ with the sheafification we obtain a map $\|K\|_\pre \times \A^1 \to \|K\|$ which, by the universal property, factors through $\|K\| \times \A^1$.
This factorization is the required weak deformation retraction of $\|K\|$ onto $\overline{\|K\|}$.
\end{proof}

\begin{proposition}
For any manifold $M$ and any presheaf $F:\man^\op\to\Sp$, the restriction map
$$\Rmap(M \times \A^\bullet,F) \to \Rmap(M \times \bDelta^\bullet,F)$$
is a weak equivalence after realization.
\end{proposition}

\begin{proof}
We can assume $F$ to be injectively fibrant, so $\Rmap$ can be replaced with $\map$.
By replacing $F$ with $\map(M,F)$ we can assume $M=\R^0$.
We will verify the conditions of \cref{prop:DuggerIsaksenSS2} and show that every square
$$\begin{tikzpicture}[descr/.style={fill=white}, baseline=(current bounding box.base)]
\matrix(m)[matrix of math nodes, row sep=2.5em, column sep=2.5em,
text height=1.5ex, text depth=0.5ex] {
\partial \Delta[n] & \Ex^\infty F(\A^\bullet) \\
\Delta[n] & \Ex^\infty \map(\bDelta^\bullet,F) \\
};
\path[->,font=\scriptsize]
	(m-1-1) edge node [auto] {$$} (m-1-2)
	(m-1-1) edge node [left] {$$} (m-2-1)
	(m-1-2) edge node [auto] {$$} (m-2-2)
	(m-2-1) edge node [below] {$$} (m-2-2);
\end{tikzpicture}$$
admits a lift making the upper and lower triangles commute up to homotopy, and such that these two homotopies are compatible on $\partial \Delta[n]$.
For a finite-dimensional simplicial compact space $K$, i.e. having a compact space of nondegenerate simplices, a map $K \to \Ex^\infty Y$ factors through some finite stage, and so it corresponds to a map $\Sd^i K \to Y$.
Then, by adjunction, the square above amounts to a map $P \to F$ where
$$P \coloneq \|\Sd^i \partial \Delta[n] \| \sqcup_{\overline{\|\Sd^i \partial \Delta[n] \|}} \overline{\|\Sd^i \Delta[n] \|}.$$
By \cref{lem:resheq}, there is a self-homotopy $H$ of $\| \Sd^i \Delta[n] \|$ such that $H_0 = \id$, $H_1$ factors through
$$\overline{\|\Sd^i \Delta[n] \|} \to  \|\Sd^i \Delta[n] \|$$
and $H$ preserves $\|\Sd^i \partial \Delta[n]\|$, $\overline{\|\Sd^i \partial \Delta[n] \|}$ and $\overline{\| \Sd^i \Delta[n] \|}$.
This defines a map $$\| \Sd^i \Delta[n]\| \to P \to F$$ and a self homotopy of $P$, giving the lift and homotopies that were needed.
\end{proof}

\subsection{Smooth maps with prescribed constancy conditions}\label{sec:constmaps}


Let $K$ be a subdivision of the standard $n$-simplex; that is, $K$ is an ordered (locally finite) simplicial complex, $|K| = \Delta^n \subset \A^n$ and every simplex of $K$ is contained (affinely) in a simplex of $\Delta^n$.
We have a map
$$j : \| K \| \rightarrow \A^n$$
that is linear on each simplex.
Note, however, that $j$ is \emph{not} induced by a simplicial map.
The map $j$ is not an inclusion but its restriction to $\overline{\| K \|}$ is, and its image is $\bDelta^n$.

\begin{proposition}\label{prop:retraction}
The inclusion of sheaves $j: \overline{\| K \|} \hookrightarrow \bDelta^n$ admits a weak deformation retraction.
More precisely, there is a map $r :\bDelta^n \to \overline{\| K \|}$ and a smooth homotopy
$$\{ h _t : \A^n \to \A^n \}_{t \in [0,1]}$$
which restricts to a homotopy $\bDelta^n \times \A^1 \to \bDelta^n$ between the identity and $jr$, and to a homotopy $\overline{\| K \|} \times \A^1 \to \overline{\| K \|}$ between the identity and $rj$.
\end{proposition}

This is a consequence of the lemma below.

\begin{lemma}\label{lem:MW}
Given a subdivision $K$ of the $n$-simplex $\Delta^n$ as in \cref{prop:retraction},
there exists a smooth homotopy
$\{ h_t : \A^n \to \A^n \}_{t \in [0,1]}$
such that
\begin{enumerate}
\item[(i)] $h_0$ is the identity,
\item[(ii)] $h_t$ maps each closed simplex $\bDelta^n \subset \overline{\|K\|}$ to itself for all $t$, and
\item[(iii)] each closed simplex $\bDelta^n \subset  \overline{\|K\|} \subset \A^n$ has a neighborhood in $\A^n$ which gets mapped to that same simplex by $h_1$.
\end{enumerate}
\end{lemma}

\begin{proof}
We use the following terminology during this proof: for a simplex $\sigma$ of $K$, a homotopy of maps $(f_t : \A^n \to \A^n)_{t \in [c,d]}$ satisfies property (iii)$_\sigma$ if $\sigma$ has a neighborhood in $\A^n$ which gets mapped to $\sigma$ by~$f_d$.

Fix some $k$ with $-1 \leq k \leq n$ and suppose per induction that we have already constructed a smooth homotopy $(h_t : \A^n \to \A^n)_{t \in [0,a]}$
for some $a < 1$ satisfying conditions (i), (ii) and (iii)$_\sigma$ for every simplex $\sigma$ of dimension at most $k$.

We want to extend this to a smooth homotopy $(h_t : \A^n \to \A^n)_{t \in [0,b]}$, where $b > a$, that satisfies conditions (i), (ii) and (iii)$_\sigma$ for every simplex $\sigma$ of dimension at most $k+1$.

For a closed $k$-simplex $\tau$, let $W_{\tau}$ be a neighborhood of $\tau$ in $\A^n$ which gets mapped to $\tau$ by $h_a$.
This exists by the inductive assumption.
We shall define a homotopy
$$(g_t : \A^n \to \A^n)_{t \in [0, b-a]}$$
where $g_0 = \id$, $g_t$ maps each simplex of $\overline{\|K\|}$ to itself for all $t$, and $g_{b-a}$ maps an appropriate subset of the interior of each $(k+1)$-simplex in $\A^n$ to that same simplex.
\emph{Appropriate} means it should be large enough so that its union with the $W_{\tau}$, over all boundary faces $\tau \subset \sigma$, contains $\sigma$, and small enough so that the various open subsets for different $(k+1)$-simplices are disjoint.
Once the homotopy $g_t$ is given, we can simply define $(h_t)_{t \in [0,b]}$ as the concatenation of $(h_t)_{t \in [0,a]}$ and $(g_{t-a} h_a)_{t \in [a,b]}$.
(In order for the concatenation to be smooth, we may arrange so that the homotopy $(h_t)_{t \in [0,a]}$ is stationary for $t$ close to $a$ and the homotopy $(g_t)$ is stationary for $t$ close to $0$.)

To describe $g_t$ we first choose, for each $(k+1)$-simplex $\sigma$, a small tubular neighborhood $U(\sigma)$ of $\int(\sigma)$ in $\A^n$
such that for each point $x \in \int(\sigma)$ and every closed simplex $\tau$ of $\|K\|$, the intersection $U(\sigma)_x \cap \tau$ is a linear cone in $\A^n$.
That is, there exist linearly independent vectors $v_1, \ldots, v_{\ell}$ such that points in $U(\sigma)_x \cap \tau$ are of the form $c_0 v_0 + c_1 v_1 + \cdots + c_{\ell} v_{\ell}$ with $c_i \geq 0$.
By shrinking if necessary, we may also assume that $U(\sigma) \cap U(\sigma')$ is empty if $\sigma$ and $\sigma'$ are distinct $(k+1)$-simplices.

Pick an open subset $V(\sigma)$ of the interior of each $(k+1)$-simplex $\sigma$, with compact closure,  whose union with $\cup_{\textup{face } \tau} W_\tau \cap \sigma$ is the closed simplex $\sigma$.
Then use the linear coordinates on the tubular neighborhood $U(\sigma)$ to obtain a map
$$\psi : U(\sigma) \to U(\sigma)$$
\emph{over} $\int(\sigma)$ satisfying the following conditions: for $x$ close to the boundary of $\sigma$, $\psi_x$ is the identity; for $x$ in $V(\sigma)$, $\psi_x(v) = 0$ for $v \in U(\sigma)_{x}$ and $|v|$ small and $\psi_x(v) = v$ for $|v|$ large.
Extend by the identity to obtain a map $g^\sigma_1 : \A^n \to \A^n$.
Linearly interpolate between the identity and $g^\sigma_1$ to get a homotopy $(g^\sigma_t)$ and concatenate the $(g^{\sigma}_t)$ for all~$\sigma$, to obtain the  homotopy $(g_t)$.
\end{proof}

\begin{proof}[Proof of \cref{prop:retraction}]
The lemma gives us a smooth homotopy $h$ on $\A^n$.
Condition (ii) implies that $h$ restricts to a homotopy $\{h_t : \overline{\|K\|} \to \overline{\|K\|} \}$.
Condition (iii) gives the required factorization of $h_1$ as
$$\bDelta^n \xrightarrow{r}  \overline{\| K \|} \xrightarrow{j} \bDelta^n$$
in the category of sheaves of sets, where the factorization of $h_1$ through $\overline{\|K\|}$ defines~$r$.
\end{proof}

\begin{remark}
\cref{lem:MW} admits a more general version which applies to arbitrary manifolds $M$ equipped with a suitable triangulation, though we will not require that level of generality.
This is claimed in Madsen--Weiss~\cite[Appendix~A.1]{Madsen.Weiss}.
\end{remark}

An inclusion of simplicial complexes $L \hookrightarrow K$ is called a \emph{relative horn inclusion} if $K$ is obtained from $L$ by attaching a simplex along a horn in $L$ (we assume that the horn is embedded in $L$).
The following lemma will be crucial to the proof of \cref{prop:CisKan}.

\begin{lemma}\label{lem:relhorninc}
Let $B$ be a subdivision of $\Delta^n$.
Given a sequence of relative horn inclusions
$$A = A_0 \hookrightarrow A_1 \hookrightarrow \cdots \hookrightarrow A_\ell = B$$
there exists a weak deformation retraction of $\overline{\|B\|}$ onto $\overline{\|A\|}$.
That is, a homotopy $H : \overline{\|B\|} \times \A^1 \to \overline{\|B\|}$
such that $H$ restricts to a homotopy $\overline{\|A\|} \times \A^1 \to \overline{\|A\|}$, $H_0 = \id$ and $H_1$ factors through $\overline{\| A \|}$.
\end{lemma}

We introduce some terminology in preparation for the proof of this lemma.
Let $0 < k \leq \ell$.
A homotopy
$$H^{(k)} : \overline{\|B\|} \times \A^1 \to \overline{\|B\|}$$
is said to have \emph{property $(\sigma_k)$} if $H^{(k)}$ restricts to a homotopy $\overline{\|A_k\|} \times \A^1 \to \overline{\|A_k\|}$, $H_0 = \id$ and $H_1$ factors through $\overline{\|A_k\|}$.

\begin{proof}
Before tackling the lemma in full generality, we prove it for the easy case of $A_0 \hookrightarrow A_1 = B$ for a single horn inclusion $\Lambda \hookrightarrow \Delta^n$.
Choose a subdivision $T$ of $\Delta^n$ which is the simplicial cone on a subdivision of $\Delta^{n-1}$.
For concreteness, we take $T$ to be the simplicial cone of $\sd \Delta^{n-1}$, the first barycentric subdivision of $\Delta^{n-1}$.
(That is, $T$ is the nerve of the category obtained by adjoining a terminal object $v$ to the poset of nondegenerate simplices of the standard $(n-1)$-simplex.)
We refer to~$T$ as the \emph{cone-subdivision} of the $n$-simplex.
Here is a picture for $n = 3$:
\begin{center}
\begin{tikzpicture}[line join = round, line cap = round]
\pgfmathsetmacro{\factor}{2/sqrt(2)}
\pgfmathsetmacro{\side}{5}
\pgfmathsetmacro{\hside}{\side/2}
\node (0) at (0,0,-1*\factor) [] {\textbullet};
\node (01) at (\hside,0.25,-1*\factor) [] {\textbullet};
\node (1) at (\side,0.5,-1*\factor) [] {\textbullet};

\node (02) at (\hside,0,0) [] {\textbullet};
\node (012) at (0.75*\side,0.15,0) [] {$c$};
\node (12) at (\side,0.3,0) [] {\textbullet};

\node (2) at (\side,0,\factor) [] {\textbullet};

\node (v) at (0.75*\side,0.7*\side,\factor) [] {$v$};

\foreach \i in {0,01,1,2,02,12,1} \draw[-] (\i) -- (012);
\foreach \i in {0,01,1,2,02,12,1,012} \draw[-] (\i) -- (v);
\draw[-] (0) -- (01);
\draw[-] (1) -- (01);
\draw[-] (0) -- (02);
\draw[-] (2) -- (02);
\draw[-] (1) -- (12);
\draw[-] (2) -- (12);
\end{tikzpicture}
\end{center}
By \cref{prop:retraction}, we have a weak deformation retraction
$$h : \bDelta^n \times \A^1 \to \bDelta^n$$
of $\bDelta^n$ onto $\overline{\|T\|}$.
Let $T'$ be the simplicial complex obtained from $T$ by discarding the vertex $c \in T$ corresponding to the top simplex in $\Delta^{n-1}$.
(To obtain a simplicial complex, we must also discard all the simplices in $T$ that have $c$ as a face.)
Then $T'$ is a subdivision of the $n$-dimensional horn and $h$ restricts to a weak deformation retraction of $\overline{\| \Lambda \|}$ onto $\overline{\|T'\|}$.

The inclusion $i : T' \hookrightarrow T$ admits a retraction $r : T \to T'$, essentially given by collapsing $c$ onto $v$.
This is a simplicial map, it is the application of the appropriate degeneracy map on each simplex of $T$.
Moreover, we can construct a homotopy on each simplex of $T$ between the identity and said degeneracy map.
This can be done by linear interpolation, for example.
Thus we obtain a deformation retraction
$$h' : \overline{\|T\|} \times \A^1 \to \overline{\|T\|}$$
of $\overline{\|T\|}$ in $\overline{\|T'\|}$.
Clearly, the composition (concatenation) of $h$ and $h'$ gives a homotopy $H$ satisfying the conditions of the lemma, i.e., having property $(\sigma_0)$.

With this special case in hand we proceed to the general one, arguing by induction.
Suppose we have constructed a homotopy $H^{(k)}$ having property $(\sigma_k)$.
We now construct a homotopy $H^{(k-1)}$ having property $(\sigma_{k-1})$ as follows.
Firstly, take a subdivision $K$ of $A_k$ (and hence a subdivision of $A_{k-1}$) that restricts to the cone triangulation on the simplex attached to $A_{k-1}$.
\cref{lem:MW} gives us a homotopy on $\bDelta^n$ that restricts to a homotopy
$$f : \overline{\|A_k\|} \times \A^1 \to \overline{\|A_k\|}$$
with $f_0 = \id$ and which factors through $\overline{\|K\|}$ at time $1$.

Now, by collapsing the cone subdivision of the attached simplex to the (subdivided) horn using the simplicial map from the case of a single horn inclusion, we obtain a homotopy
$$g : \overline{\|K\|} \times \A^1 \to \overline{\|K\|} \subset \overline{\|A_k\|}$$
with $g_0 = \id$ and which factors through $\overline{\|A_{k-1}\|}$ at time $1$.
Compose (concatenate) the homotopies $f$ and $g$ to obtain a homotopy $h$ on $\| A_k \|$.
Then define $H^{(k-1)}$ to be composition of $H^{(k)}$ and $h$.
For this composition to be smooth, we emphasize that it is important to first apply \cref{lem:MW} to the whole of $\overline{\|A_k\|}$, not just the on the simplex that we are collapsing.
This completes the induction.
\end{proof}

\begin{corollary}\label{cor:relhorninc}
The extended-simplices version of \cref{lem:relhorninc} holds.
Namely, under the conditions of that lemma, $\| B \|$ weak deformation retracts onto $\| A \|$.
That is, there exists a homotopy $G : \| B \| \times \A^1 \to \| B \|$
such that $G$ restricts to a homotopy $\| A \| \times \A^1 \to \| A \|$, $G_0 = \id$ and $G_1$ factors through $\| A \|$.
\end{corollary}

\begin{proof}
Starting from the homotopy $H : \overline{\| B \|} \times \A^1 \to \overline{\| B \|}$ of \cref{lem:relhorninc}, we obtain a homotopy
$$\widetilde{H} : \| B \| \times \A^1 \xrightarrow{\lambda} \overline{\| B \|} \times \A^1 \xrightarrow{H} \overline{\| B \|} \xrightarrow{i} \| B \|$$
where $i$ is the inclusion and $\lambda$ is the map constructed in \cref{lem:resheq}.
It is clear that $\widetilde{H}$ restricts to a homotopy on $\| A \|$, $\widetilde{H}_0 = \lambda i$ and $\widetilde{H}_1$ factors through $\| A\|$.
Now define $G$ as the concatenation of the homotopy on $\| B \|$ from \cref{lem:resheq} (between the identity and $ \lambda i$) with the homotopy $\widetilde{H}$.
\end{proof}

\subsection{Proof of \cref{prop:CisKan,prop:EPisKan,prop:SnisKan}}\label{sec:EPisKanassrelhorn}

For concreteness, we assume $F$ to be injectively fibrant (by replacing it if necessary) so that \cref{prop:inj} applies.

\begin{proof}[Proof of \cref{prop:CisKan}]
Let $I:  A \hookrightarrow B$ denote the map
$$h_i:\Sd^i (\Lambda_k[n])  \hookrightarrow \Sd^i (\Delta[n])$$
with $i \geq 0$, $n \geq 1$, and $0 \leq k \leq n$, as in \cref{defn:wKan}.
We write $\iota=\|I\|_\pre$.
By adjunction of \cref{adjunction.pre}, it suffices to find weak liftings
$$\begin{tikzpicture}[descr/.style={fill=white}, baseline=(current bounding box.base)]
\matrix(m)[matrix of math nodes, row sep=2.5em, column sep=2.5em,
text height=1.5ex, text depth=0.25ex] {
\| A \|_\pre & F \\
\| B \|_\pre & \\
};
\path[->,font=\scriptsize] (m-1-1) edge node [auto] {$f$} (m-1-2);
\path[->,font=\scriptsize] (m-1-1) edge node [left] {$\iota$} (m-2-1);
\path[->,dashed,font=\scriptsize] (m-2-1) edge node [auto] {$\widetilde{\alpha}$} (m-1-2);
\end{tikzpicture}$$
together with a homotopy $H : \| A \|_\pre \times \| \Delta^1 \|_\pre \to F$ between $\widetilde{\alpha} \iota$ and $f$.
Indeed, this gives us the required homotopy $A \times \Delta[1] \to F(\A^\bullet)$ by precomposing $H$ with $\| A \times \Delta^1 \|_\pre \to  \| A \|_\pre \times \| \Delta^1 \|_\pre$ and applying the adjunction again.
It also suffices to solve the above problem with $\|{-}\|_\pre$ replaced by $\|{-}\|$ everywhere.
This is allowed since $F$ is an $\infty$-sheaf and the map $\|{-}\|_\pre\to\|{-}\|$ is by definition a sheafification.

The strategy is to find a homotopy retraction of $\iota$, i.e., a map $r : \| B \|\to \| A \|$ together with a homotopy $ \| A \| \times \| \Delta^1 \| \to  \| A \|$ between $r \iota$ and the identity.
In fact, we construct this homotopy on $\| A \|$ as the restriction of a homotopy on $\| B \|$.
(We will need this stronger statement in the proof of \cref{prop:EPisKan}.)
This is achieved by a direct application of \cref{cor:relhorninc}.
More precisely, we choose a sequence of relative horn inclusions from $A = \sd^i \Lambda^n_k$ to $B = \sd^i \Delta^n$, for each $i \geq 0$, $n > 0$ and $0 \leq k \leq n$, and apply \cref{cor:relhorninc}.
(A proof that such a sequence exists can be found in Moss~\cite[Proposition~19]{Moss}.)
That the resulting homotopy preserves the $0$-simplices of each horn follows directly from its construction.
\end{proof}

\begin{proof}[Proof of \cref{prop:EPisKan}]
We keep the notation $I : A \hookrightarrow B$ for the map~$h_i$ as in \cref{defn:wKan}.
As in the proof of \cref{prop:CisKan}, we use the adjunction of \cref{adjunction.pre} and the natural transformation $\|{-}\|_\pre\to\|{-}\|$ to reduce
the problem to constructing certain maps of sheaves of sets.

Let $P$ denote the pushout of sheaves of sets
$$\| A \| \times \A^1 \sqcup_{\|A \| \times \partial \A^1}  \| B \|\times \partial \A^1.$$
To verify the weak RLP with respect to $\|I\|$, it suffices to prove that for any map $\alpha:P\to F$ there is a dashed map $\widetilde{\alpha}$ as in
$$\begin{tikzpicture}[descr/.style={fill=white}, baseline=(current bounding box.base)]
\matrix(m)[matrix of math nodes, row sep=2.5em, column sep=2.5em,
text height=1.5ex, text depth=0.5ex] {
P & F \\
\| B \| \times \A^1 & \\
};
\path[->,font=\scriptsize] (m-1-1) edge node [auto] {$\alpha$} (m-1-2);
\path[->,font=\scriptsize] (m-1-1) edge node [left] {$\iota$} (m-2-1);
\path[->,dashed,font=\scriptsize] (m-2-1) edge node [below] {$\widetilde{\alpha}$} (m-1-2);
\end{tikzpicture}$$
making the diagram commute up to a homotopy $H : P \times \A^1 \to F$ from $\alpha$ to $\widetilde{\alpha} \iota$, which is \emph{fixed} on $\| B \|\times\partial \A^1$ pointwise.
Being fixed pointwise means that the restriction of $H$ to $(\| B \|\times\partial \A^1) \times \A^1 $ factors as the projection to $\| B \|\times\partial \A^1$ followed by $\alpha$.
The result now follows from the lemma below.
\end{proof}

\begin{lemma}\label{lem:410}
Suppose $I : A \hookrightarrow B$ is the map~$h_i$ as in \cref{defn:wKan} and $P$ is the pushout of sheaves of sets $\| A \| \times \A^1 \sqcup_{\|A \| \times \partial \A^1}  \| B \|\times \partial \A^1$.
Then $\| B \|\times \A^1$ weak deformation retracts onto $P$, relative to $\| B \|\times \partial \A^1$.
\end{lemma}
\begin{proof}
We need to show that there is a homotopy $H : (\| B \|\times \A^1) \times \A^1 \to \| B \|\times \A^1$ such that
\begin{enumerate}
\item[(1)] $H$ restricts to a homotopy $P \times \A^1 \to P$ which fixes $\| B \|\times \partial \A^1$ pointwise
\item[(2)] $H_0$ is the identity and $H_1$ factors through $P$.
\end{enumerate}

Choose a bump function $c: \A^1\to [0, 1] \subset \A^1$ with $c\equiv0$ in an open neighborhood of~$(-\infty,0] \cup [1, \infty)$, and $c\equiv 1$ in a open neighborhood $J$ of $1/2$, $c(t)$ increasing for $t \leq 1/2$ and decreasing for $t \geq 1/2$.
Then choose a function $f: \A^1\to \A^1$ with $f(t)\equiv 0$ when $c\ne 1$ and $t\le 1/2$ and $f(t)\equiv 1$ when $c\ne 1$ and $t>1/2$.
Let $H_B: \|B\|\times \A^1\to \|B\|$ be the weak deformation retraction constructed in the proof of \cref{prop:CisKan}.
Define maps $R_1, R_2 : (\|B\|\times \A^1)\times \A^1 \to \|B\|\times \A^1$ as
$$R_1(x,t,s) = (H_B(x,s\cdot c(t)),t)$$
and
$$R_2(x,t,s) = (x,s\cdot f(t)+(1-s)t).$$
for $x \in {\| B\|}$, and $t, s \in \A^1$.
Then define $H$ as
$$(x,t,s)\mapsto \begin{cases}
R_1(x,t,2s) & \text{for } s\le 1/2 \cr
R_2(R_1(x,t,1),2s-1) & \text{for } s>1/2 \cr
\end{cases}$$
Now $R_1$ and $R_2$ separately satisfy condition~(1), so $H$ does as well.
As for property~(2), we have that $H(x,t,1) = (x', t')$ where $x' = H_B(x, c(t))$ and $t' = f(t)$.
If $t \in J$ then $c = 1$ and so $x' \in \| A \|$, and if $t \notin J$ then $f(t) \in \partial \A^1$.
Therefore $H_1$ factors through $P$.
\end{proof}

\begin{proof}[Proof of \cref{prop:SnisKan}]
The proof is of the same sort as that of \cref{prop:EPisKan},
using the adjunction of \cref{adjunction.pre} and the natural transformation $\|{-}\|_\pre\to\|{-}\|$ to reduce
the problem to constructing certain maps of sheaves of sets.
Let $I : A \hookrightarrow B$ be the map~$h_i$ as in \cref{defn:wKan}.
Consider the pushout
$$P \coloneq S^n\times \| A\| \sqcup_{*\times \|A\|} *\times \|B\|.$$
The same manipulations as before show that to verify the weak RLP with respect to $\|I\|$ it suffices to prove that for any map $\alpha:P\to F$ there is a map~$\widetilde{\alpha}$ as in
$$\begin{tikzpicture}[descr/.style={fill=white}, baseline=(current bounding box.base)]
\matrix(m)[matrix of math nodes, row sep=2.5em, column sep=2.5em,
text height=1.5ex, text depth=0.5ex] {
P & F \\
S^n\times \|B\| & \\
};
\path[->,font=\scriptsize] (m-1-1) edge node [auto] {$\alpha$} (m-1-2);
\path[->,font=\scriptsize] (m-1-1) edge node [left] {$\iota$} (m-2-1);
\path[->,dashed,font=\scriptsize] (m-2-1) edge node [below] {$\widetilde{\alpha}$} (m-1-2);
\end{tikzpicture}$$
making the diagram commute up to a homotopy $P\times \A^1\to F$ from $\alpha$ to $\tilde\alpha\iota$ which is fixed on $*\times \|B\|$ pointwise.
The result now follows from a lemma analogous to \cref{lem:410}, below.
\end{proof}

\begin{lemma}
Suppose $I : A \hookrightarrow B$ is the map~$h_i$ as in \cref{defn:wKan} and $P$ is the pushout of sheaves of sets
$S^n\times \| A\| \sqcup_{*\times \|A\|} *\times \|B\|$.
Then $\| B \times S^n \|$ weak deformation retracts onto $P$, relative to $\| B\| \times *$.
\end{lemma}

\begin{proof}
As before, we need to show that there is a homotopy $H : (\| B \|\times S^n) \times \A^1 \to \| B \|\times S^n$ such that
\begin{enumerate}
\item $H$ restricts to a homotopy $P \times \A^1 \to P$ which fixes $\| B \|\times *$ pointwise.
\item $H_0$ is the identity and $H_1$ factors through $P$.
\end{enumerate}

Let $*\in D_\epsilon\subset D_\delta\subset S^n$ be open disk neighborhoods of radii $\epsilon$ and $\delta$ with $\epsilon<\delta$.
Choose a smooth function $c: S^n\to [0,1]\subset \A^1$ such that $c|_{S^n\setminus D_\epsilon}\equiv 0$ and $c(*)=1$.
For example, we can choose $c$ to be a bump function with support in $D_\epsilon $ that is $1$ at the basepoint.

Also choose a homotopy $h: \A^1\times S^n\to S^n$ such that $h(1,-)|_{D_\epsilon}$ is constant to $*\in S^n$ and $h(t,-)|_{S^n\setminus D_\delta}=\id$.
In words, this homotopy collapses a neighborhood of the basepoint to itself.

Let $H_B: \|B\|\times \A^1\to \|B\|$ be the map constructed in the proof of \cref{prop:CisKan}.
We construct $H$ as the composition of two homotopies.
Define
$$R_1: S^n\times \|B\|\times \A^1\to S^n\times \|B\| \qquad R_1(z,x,t)=(z,H_B(x,t(1-c(z))))$$
and
$$R_2: S^n\times \|B\|\times \A^1\to S^n\times \|B\| \qquad R_2(z,x,t)=(h(t,z),x).$$
The first homotopy collapses $S^n \times \|B\|$ onto $S^n \times \|A\|$ outside a neighborhood $N\subset D_\epsilon$ of the basepoint.
The second collapses $D_\epsilon\times \|B\|$ to $*\times \|B\|$.
The composition of these homotopies satisfies the claimed properties.
\end{proof}

\section{The shape functor preserves the $\infty$-sheaf property}
\label{sec:sheafprop}

In this section we assemble the previous results to prove \cref{sheafthm}.
Our approach uses the following characterization of $\infty$-sheaves.

\begin{theorem} \label{thm:sheafcond}
A presheaf $F:\man^\op\to\Sp$ is an $\infty$-sheaf if and only if $F(\varnothing) \simeq *$ and
\begin{enumerate}
\item[(1)] for all manifolds $M$ and open covers of $M$ with two elements $\{U,V\}$,
the commutative square
\begin{equation*}
\begin{tikzpicture}[descr/.style={fill=white}, baseline=(current bounding box.base)]
\matrix(m)[matrix of math nodes, row sep=2.5em, column sep=2.5em,
text height=1.5ex, text depth=0.25ex] {
F(M) & F(V) \\
F(U) & F(U \cap V) \\
};
\path[->,font=\scriptsize] (m-1-1) edge node [auto] {} (m-1-2);
\path[->,font=\scriptsize] (m-2-1) edge node [auto] {} (m-2-2);
\path[->,font=\scriptsize] (m-1-1) edge node [left] {} (m-2-1);
\path[->,font=\scriptsize] (m-1-2) edge node [auto] {} (m-2-2);
\end{tikzpicture}
\end{equation*}
is a homotopy pullback square; and
\item[(2)] If $M$ is a (possibly infinite) disjoint union of submanifolds $U_i$ then $F(M) \to \prod^h_{i} F(U_i)$ is a weak equivalence.
\end{enumerate}
\label{thm:weiss} \end{theorem}

This is probably well-known and is similar to a special case of Weiss~\cite[Theorem~5.2]{Weiss.Embeddings} and Boavida–Weiss \cite[Theorem~7.2]{Boavida.Weiss}.
For completeness, we provide a proof below.
In preparation, we record the following:
\begin{lemma}\label{lem:holims}
A presheaf $F : \man^\op\to\Sp$ is an $\infty$-sheaf if and only if for every open cover $\{U_i\}_{i \in I}$, the canonical map
$$F(M) \to \holim_{S \subset I} F(U_{S})$$
is a weak equivalence, where the homotopy limit ranges over all \emph{finite, nonempty subsets} $S \subset I$ and $U_S$ is notation for $\bigcap_{i \in S} U_i$.
\end{lemma}
\begin{proof}
For the duration of this proof we will write $\underline n$ for the set $\{0,\dots, n\}$ to distinguish it from the total ordered set $[n] \coloneq \{0 \leq \dots \leq n\}$.

We need to show that 
\begin{equation}\label{eq:holims}
\holim_{[n] \in \Delta} \prod_{i_0, \ldots, i_n \in I} F(U_{i_0} \cap \cdots \cap U_{i_n})\simeq \holim_{S \subset I} F(U_{S}) \; .
\end{equation}

First suppose that $I$ is finite and pick a total order on $I$. We can then replace the homotopy limit on the left by replacing the product over sequences $i_0, \dots, i_n$ with the product over ordered sequences ${i_0 \leq \dots \leq i_n}$. To see this, view the left-hand side of the display as a homotopy limit of a (covariant) functor from $D$ to spaces where $D$ is the category whose objects are sequences $i_0, \dots, i_n$ and morphisms are induced by $\Delta$ and correspond to merging repeated elements or adding new ones.  More precisely, an object is a map $i : \underline n \to I$ and a morphism from $t : \underline m \to I$ to $i : \underline n \to I$ is an order preserving map $\theta : \underline m \to \underline n$ such that $t = \theta i$.

Let $D'$ be the full subcategory of $D$ consisting of \emph{ordered} sequences $i_0 \leq \dots \leq i_n$, i.e. functors $[n] \to I$.  There is a canonical functor $o : D \to D'$ that orders each sequence. Namely, for a sequence $i_0, \dots, i_n$, i.e. a map $i : \underline n \to I$, precompose with the unique bijection $f :\underline n \to \underline n$ such that $if$ is order-preserving and the restriction of $f$ to $(if)^{-1}(j)$ is order-preserving for each $j \in I$. Given a morphism $\theta$ in $D$ from $t : \underline m \to I$ to $i : \underline n \to I$, its image under $o$ is the morphism from $o(t) = gt$ to $o(i) = if$ in $D^\prime$ given by $f^{-1} \theta g : [m] \to [n]$ (this is indeed order-preserving since $g$ reverses the order of two elements if and only if $f$ does).

We claim that this functor $o$, sending $i$ to $if$, is homotopy initial. That is, for each ordered sequence $[\underline i] \coloneq (i_0 \leq \dots \leq i_n)$, the comma category $o/[\underline i]$ is contractible. To prove this, we show that the identity map on the classifying space of $o/[\underline i]$ is null-homotopic by considering functors $A, \textup{const} : o/[\underline i] \to o/[\underline i]$. The functor $A$ sends an object in $o/[\underline i]$, i.e. a morphism $o(a_0, \dots, a_k) \to [\underline i]$ in $D^\prime$, to the object $o(i_0, \dots, i_n, a_0, \dots, a_k) \to [\underline i]$ which merges all repetitions. The functor \textup{const} sends all objects to the identity $o(i_0, \dots, i_n) \to [\underline i]$. For each object $o(a_0, \dots, a_k) \to [\underline i]$, there are morphisms in $o/[\underline i]$
\[
o(a_0, \dots, a_k) \to o(i_0, \dots, i_n, a_0, \dots, a_k) \gets o(i_0, \dots, i_n)
\]
induced by faces in $\Delta$, i.e. adding new elements. These define natural transformations $\id \Rightarrow A \Leftarrow \textup{const}$, proving the claim.

Secondly, given an ordered sequence we can forget its ordering and view it as a finite subset of $I$.
This construction defines a functor~$\eta$ from $D'$ to the category of nonempty subsets of $I$.
{Then, for each $S \subset I$, the comma category $\eta/S$ is contractible since it is the category of simplices of the nerve of $S$ (viewed as a poset with respect to the total ordering induced from the inclusion $S \subset I$). In other words, $\eta$ is homotopy initial.}

We have shown that, for a finite set $I$, the map induced by $\eta o$ gives a weak equivalence (\ref{eq:holims}). It is clear that this is natural with respect to inclusions of finite sets $I^\prime \subset I$. So to prove (\ref{eq:holims}) for a general indexing set $I$ we can reduce to the finite case by taking the homotopy limit of weak equivalences (\ref{eq:holims}) over all finite subsets of $I$. In more detail, we have a commutative square
\begin{equation*}
\begin{tikzpicture}[descr/.style={fill=white}, baseline=(current bounding box.base)]
\matrix(m)[matrix of math nodes, row sep=2.5em, column sep=2.5em,
text height=1.5ex, text depth=0.25ex] {
\holim_{S \subset I} F(U_S) & \holim_{J \subset I} \holim_{S \subset J} F(U_S) \\
\holim_{i : \underline n \to I} F(U_{i_0} \cap \cdots \cap U_{i_n}) & \holim_{J \subset I} \holim_{i : \underline n \to J} F(U_{i_0} \cap \cdots \cap U_{i_n}) \\};
\path[->,font=\scriptsize] (m-1-1) edge node [auto] {} (m-1-2);
\path[->,font=\scriptsize] (m-2-1) edge node [auto] {} (m-2-2);
\path[->,font=\scriptsize] (m-1-1) edge node [auto] {} (m-2-1);
\path[->,font=\scriptsize] (m-1-2) edge node [auto] {} (m-2-2);
\end{tikzpicture}
\end{equation*}
where the first homotopy limits in the right column run over finite subsets $J \subset I$.
We have shown that the right-hand map is a weak equivalence.
The horizontal maps are also weak equivalences which we can see by identifying the double homotopy limits with a single homotopy limit over a Grothendieck construction (Thomason's homotopy colimit theorem).
We will explain this briefly for the lower map, the upper one is similar and almost identical to an argument in the proof of \cref{thm:sheafcond}.
In that case, the double homotopy limit is identified with the homotopy limit over the category whose objects are pairs $\underline n \to J \subset I$ with $J$ finite, and morphisms are maps of such.
The forgetful functor from this category to $D$, sending a pair $\underline n \to J \subset I$ to its composite $\underline n \to I$, is homotopy initial since the relevant comma categories have a terminal object.
\end{proof}

\begin{proof}[Proof of \cref{thm:sheafcond}]

If $F$ is an $\infty$-sheaf then conditions (1) and (2) are immediate in view of \cref{lem:holims}.

To show the converse, suppose first that $M$ is a compact manifold and take an open cover $\{U_i\}_{i \in I}$ of $M$.
For every finite subcover $\{U_j\}_{j \in J}$ with $J \subset I$ of $\{U_i\}_{i \in I}$, the homotopy limit
$$\holim_{S \subset J} F(U_S)$$
is indexed over a finite category (a cube) and so it is equivalent to an iterated homotopy pullback.
Condition~(1) applied inductively shows that this iterated homotopy pullback is $F(M)$.
Then consider the square
\begin{equation*}
\begin{tikzpicture}[descr/.style={fill=white}, baseline=(current bounding box.base)]
\matrix(m)[matrix of math nodes, row sep=2.5em, column sep=2.5em,
text height=1.5ex, text depth=0.25ex] {
F(M) & \holim_{J \subset I} F(M) \\
\holim_{S \subset I} F(U_S) & \holim_{J \subset I} \holim_{S \subset J} F(U_S) \\
};
\path[->,font=\scriptsize] (m-1-1) edge node [auto] {$\simeq$} (m-1-2);
\path[->,font=\scriptsize] (m-2-1) edge node [auto] {$\simeq$} (m-2-2);
\path[->,font=\scriptsize] (m-1-1) edge node [left] {} (m-2-1);
\path[->,font=\scriptsize] (m-1-2) edge node [auto] {$\simeq$} (m-2-2);
\end{tikzpicture}
\end{equation*}
where the outer homotopy limits in the right column are indexed by finite refinements, i.e., finite subsets $J \subset I$ such that $\{U_j\}_{j \in J}$ is still a cover.
The right-hand map is a weak equivalence by the observation just made.
The poset of finite refinements is filtered, and hence contractible, and so the top horizontal map in the square is also a weak equivalence.
As for the lower horizontal map one can, by Thomason's homotopy colimit theorem, express the double homotopy limit as a homotopy limit over the category---call it $P$---whose objects are pairs $S \subset J \subset I$ and morphisms are inclusions of such.
The forgetful map $\eta$ from $P$ to the poset of finite, nonempty subsets $S \subset I$ is homotopy initial since, for each $S \subset I$, the overcategory $\eta/S$ is the filtered poset of all refinements $J$ containing $S$.

To prove the noncompact case, we will assume for the moment that $F$ satisfies the following condition, which implies condition~(2):
\begin{enumerate}
\item[($2'$)] for any manifold $M$ and an open cover $\{V_i\}_{i\ge0}$ by a nested sequence of open sets with $\overline{V_i} \subset V_{i+1}$, the canonical map
$$F(M)=F \left(\bigcup_i V_i\right)\to \holim_i F(V_i)$$ is a weak homotopy equivalence.
\end{enumerate}
Assuming $F$ satisfies (1) and ($2'$), we will now prove that it satisfies the $\infty$-sheaf condition for any noncompact manifold.
So let $M$ be a noncompact manifold, and take an exhaustion of $M$ by interiors of compact manifolds $V_0 \subset V_1 \subset \cdots$ with $M = \bigcup V_i$.
Such an exhaustion can be obtained by picking a smooth proper map $f : M \to \R$ and setting $V_i$ to be the interior of $f^{-1}((-\infty,i])$.
Then, for an open cover $\{U_i \to M\}_{i \in I}$,
\begin{equation}\label{eq:bla}
\holim_{S \subset I} F(U_S) \simeq \holim_{S \subset I} \holim_{j \geq 0} F(V_j \cap U_S)
\end{equation}
by ($2'$) applied to the covers $\{U_S \cap V_j \to U_S\}_{j}$ for each $S$.
Now commute the homotopy limits and use that the cover $\{V_j \cap U_i \to V_j\}_i$ (for a fixed $j$) has a finite subcover to conclude using (1) that (\ref{eq:bla}) is weakly equivalent to
$$\holim_{j \geq 0} F(V_j).$$
By invoking ($2'$) again, this homotopy limit is weakly equivalent to $F(M)$.

We are left to show that conditions~(1) and (2) jointly imply condition~($2'$).
Suppose $\{V_i\}$ is an open cover as in~($2'$).
Let $W_0$ be the disjoint union of $V_{k+2} \backslash \overline{V_k}$, taken over $k$ even, and let $W_1$ be the disjoint union of $V_{k+2} \backslash \overline{V_k}$, taken over $k$ odd.
Then $W_0$ and $W_1$ form an open cover of $M$ so, by (1), we have an equivalence
$$F(M) \to \holim (F(W_0) \to F(W_0 \cap W_1) \gets F(W_1)).$$
By ($2$), the target is equivalent to
$$\holim ( \holim_{i} F(W_0 \cap V_i) \to \holim_{i}F(W_0 \cap W_1 \cap V_i) \gets \holim_{i}F(W_1 \cap V_i)),$$
which, by commuting homotopy limits and using (1), is equivalent to $\holim_{i} F(V_i)$.
\end{proof}

\begin{remark}\label{rem:sheafcond}
In a previous iteration of this paper, \cref{thm:sheafcond} had a stronger condition (2).
The referee kindly pointed out to us that our proof implied the new (weaker) statement and suggested the simple argument of the last paragraph of the proof.
\end{remark}

\begin{remark}
The same proof works for topological manifolds.
The main observation for the noncompact case is that there exists a proper map $M \to \R$ for $M$ a topological manifold (the requirement is that partitions of unity exist).
For generalizations of these statements, see Pavlov~\cite{Pavlov.BG}.
\end{remark}

We shall tackle properties (1) and (2) for $\sB F$ separately below.
We call them the \emph{finite} and \emph{noncompact} cases, respectively.

\subsection{The finite case}
\begin{theorem}\label{thm:compact}
Let $F:\man^\op\to\Sp$ be an $\infty$-sheaf and $M$ a smooth manifold with $U$ and $V$ two open subsets of~$M$ such that $U \cup V = M$.
Then the commutative square
\begin{equation*}
\begin{tikzpicture}[descr/.style={fill=white}, baseline=(current bounding box.base)]
\matrix(m)[matrix of math nodes, row sep=2.5em, column sep=2.5em,
text height=1.5ex, text depth=0.25ex] {
\sB F(M) & \sB F(V) \\
\sB F(U) & \sB F(U \cap V) \\
};
\path[->,font=\scriptsize] (m-1-1) edge node [auto] {} (m-1-2);
\path[->,font=\scriptsize] (m-2-1) edge node [auto] {} (m-2-2);
\path[->,font=\scriptsize] (m-1-1) edge node [left] {} (m-2-1);
\path[->,font=\scriptsize] (m-1-2) edge node [auto] {} (m-2-2);
\end{tikzpicture}
\end{equation*}
is homotopy cartesian.
\end{theorem}

\begin{proof}
The (homotopy) pullback
$$(\sB F(U) \times \sB F(V))\times_{\map(\partial \Delta^1,\sB F(U \cap V))} \map(\Delta^1,\sB F(U \cap V))$$
is identified with
\begin{equation}\label{eq:hpb}
(\sB F(U) \times \sB F(V)) \times^h_{\sB F(U \cap V \times \partial \A^1)} \sB F(U\cap V \times \A^1)
\end{equation}
since $\sB F$ is concordance invariant.
By \cref{prop:EPisKan}, we may commute the homotopy pullback with geometric realization, and thus (\ref{eq:hpb}) is identified with the geometric realization of the simplicial space
\begin{equation}\label{eq:hpb2}
(F(U \times \A^\bullet) \times F(V \times \A^\bullet)) \times^h_{F(U \cap V \times \partial \A^1 \times \A^\bullet)}
F(U \cap V \times \A^1 \times \A^\bullet).
\end{equation}
To prove that the map from $F(M \times \A^{\bullet})$ to (\ref{eq:hpb2}) is a weak equivalence after realization we first refine the cover in a convenient way using a partition of unity subordinate to $\{U,V\}$.
So let $f_U: M\to [0,1]$ and $f_V: M\to [0,1]$ with $f_U+f_V\equiv 1$, and ${\rm supp}(f_U)\subset U$ and ${\rm supp}(f_V)\subset V$.
Take $U'=f_U^{-1}(2/3,1]$ and $V'=f_V^{-1}(2/3,1]$.
Notice that $U'\cap V'=\emptyset$, and $\{U',V',U\cap V\}$ covers~$M$.
Let $c: \A^1 \to \A^1$ be a cutoff function with $c|_{(-\infty,1/3)}\equiv 0$ and $c|_{(2/3,\infty)}\equiv 1$, and define~$f\coloneq c\circ f_V|_{U\cap V}$, so that $f:U\cap V\to\A^1$.

Rearrange (\ref{eq:hpb2}) as an iterated homotopy pullback and consider the maps:
$$\begin{tikzpicture}[node distance=3.5cm,auto]
\node (C) [node distance = 6.5cm] {$ F(U \times \A^\bullet)\times^h_{F(U\cap V \times \A^\bullet)} F(U\cap V\times \A^1 \times \A^\bullet)\times^h_{F(U\cap V \times \A^\bullet)} F(V \times \A^\bullet)$};
\node (B) [node distance= 1.5cm, below of=C] {$F(U' \times \A^\bullet)\times^h_{F(U'\cap V \times \A^\bullet)} F(U \cap V \times \A^1 \times \A^\bullet)\times^h_{F(U\cap V' \times \A^\bullet)}F(V' \times \A^\bullet)$};
\node (BB) [node distance=4.5cm, left of=B] {$$};
\node (D) [node distance =1.5cm, below of =B] {$F(U' \times \A^\bullet)\times^h_{F(U'\cap V \times \A^\bullet)} F(U \cap V \times \A^\bullet)\times^h_{F(U \cap V' \times \A^\bullet)}F(V' \times \A^\bullet).$};
\node (DD) [node distance = 4.2cm, left of=D] {$$};
\draw[->] (C) to node {${\rm res}$} (B);
\draw[->] (D) to node {${\rm pr}^*$} (B);
\draw[->, bend left] (B) to node {$f^*$} (D);
\end{tikzpicture}$$
The restriction map from $F(M \times \A^\bullet)$ to the last space is a weak equivalence since $\{U', U \cap V, V'\}$ covers $M$ and $F$ is an $\infty$-sheaf.
Similarly, the map ${\rm res}$ is a levelwise weak equivalence since $\{U',U\cap V\}$ covers~$U$ and $\{V',U\cap V\}$ covers~$V$.
The arrow ${\rm pr}^*$ is induced by the projection ${\rm pr}: (U\cap V)\times \A^1 \to U\cap V$.
We obtain a map $U\cap V\to (U\cap V)\times \A^1$ from~$f: U\cap V\to \A^1$ defined in the previous paragraph.
By construction, $f|_{U'\cap V}=0$ and $f|_{U\cap V'}=1$, which is precisely the compatibility condition required to extend to a map on sections in the fibered product which we denote by~$f^*$.
Notice that since ${\rm pr}\circ (\id_{U\cap V},f)=\id_{U\cap V}$, we have $(\id_{U\cap V},f)^*\circ {\rm pr}^*=\id$.
It remains to show that ${\rm pr}^*\circ (\id_{U\cap V},f)^*$ is homotopic to the identity.

We consider the interpolation $h: ((U\cap V)\times \A^1) \times \A^1\to \A^1$ between~$f\circ {\rm pr}$ and the projection map $q:(U\cap V)\times \A^1\to \A^1$, given by
$$h(t)=(1-t)\cdot q+t\cdot (f\circ {\rm pr})$$
and extend it to a smooth homotopy $$H = (\id_{U\cap V}, h) : ((U\cap V)\times \A^1) \times \A^1 \to ((U\cap V)\times \A^1).$$
Since $F(-\times \A^\bullet)$ sends smooth homotopies to simplicial homotopies (\cref{prop:enriched}), and the map~$H$ fixes $(U\cap V) \times \partial \A^1$ pointwise,
the map~$H$ induces the required simplicial homotopy from ${\rm pr}^*\circ (\id_{U\cap V},f)^*$ to $\id$.
\end{proof}

\begin{corollary}\label{cor:compact}
Let $F$ be a presheaf on $\man$ which satisfies the $\infty$-sheaf condition with respect to finite covers.
Then the evaluation map
$$\sB F(M) \rightarrow \map(\Sing M, \tB F)$$
is a natural weak equivalence of spaces for every compact manifold $M$.
\end{corollary}
\begin{proof}
This follows from \cref{thm:compact} and the proof of \cref{prop:Dugger} applied to a finite good open cover of $M$.
\end{proof}

The beginning of the proof of \cref{thm:compact} has the following obvious generalization, which is just a consequence of \cref{prop:EPisKan}.

\begin{definition}
For a diagram $F \to G \gets H$ of $\infty$-sheaves, define the \emph{geometric} homotopy pullback to be the $\infty$-sheaf whose value at a manifold $M$ is the homotopy pullback of the diagram
$$\begin{tikzpicture}[descr/.style={fill=white}, baseline=(current bounding box.base)]
\matrix(m)[matrix of math nodes, row sep=2.5em, column sep=2.5em,
text height=1.5ex, text depth=0.25ex] {
& G(M \times \A^1) \\
F(M) \times H(M) & G(M) \times G(M) \\
};
\path[->,font=\scriptsize]
	(m-2-1) edge node [auto] {} (m-2-2)
	(m-1-2) edge node [auto] {$\textup{endpoints}$} (m-2-2);
\end{tikzpicture}$$
\end{definition}
Then, by \cref{prop:EPisKan}, the classifying space functor $\sB$ sends geometric homotopy pullbacks of $\infty$-sheaves to homotopy pullbacks of spaces.

\subsection{The noncompact case}

\begin{theorem}\label{thm:noncompact}
Let $\{U_i\}_{i \geq 0}$ be a collection of manifolds (possibly noncompact).
For an $\infty$-sheaf $$F:\man^\op\to\Sp,$$ the natural map
$$\sB F\Bigl( \bigsqcup_{i} U_i\Bigr) \to \prod^h_i \sB F(U_i)$$
is a weak equivalence.
\end{theorem}

We deduce this from the lemma below, by setting $F_i = F(U_i \times -)$.
\begin{lemma}\label{cor:infinite-prods}
Let $\{F_i\}_{i \in I}$ be a collection of $\infty$-sheaves indexed over a possibly infinite set $I$.
Then the map
\begin{equation}\label{eq:inftyprods}
\biggl|\prod^h_{i} F_i(\A^\bullet)\biggr| \to \prod_i^h |F_i(\A^\bullet)|
\end{equation}
is a weak equivalence of spaces.
In other words, the functor~$\tB$ preserves small homotopy products. 
\end{lemma}

We use the symbol $\prod^h$ for the homotopy product, i.e., the derived functor of the product.
This has a different meaning in simplicial spaces (with degreewise weak equivalences) and simplicial sets (with the usual weak equivalences).
In the simplicial space case, it means: replace each factor by a degreewise fibrant simplicial space and then compute the product;
in the simplicial set case, it means: replace each factor by a Kan complex and then compute the product.
The homotopy product (of spaces or of simplicial spaces) agrees with the nonderived product when the indexing set is finite.
In general, they do not agree when the set is infinite but \cref{cor:infinite-prods} says they agree for the concordance resolution.

\begin{proof}[Proof of \cref{cor:infinite-prods}]
The following elegant argument was suggested to us by a referee.
Without loss of generality, we can assume $F_i$ to be injectively fibrant, by performing an injective fibrant replacement if necessary.
Thus, $F_i$ is valued in Kan complexes and the homotopy product $\prod^h_i F_i$ can be computed as the ordinary product $\prod_i F_i$.
It then suffices to show that the map (\ref{eq:inftyprods}) induces an isomorphism on homotopy groups for all degrees and base-points.

We have
$$\pi_n\Bigl( \prod_{i} | F_i(\A^\bullet) |, x\Bigr) \cong \prod_{i} \pi_n (|F_i(\A^\bullet)|, x_i) \cong \prod_{i} \pi_0 F_i(S^n)_{x_i} / {\sim}$$
where $F_i(S^n)_{x_i}$ is the fiber $F_i(S^n) \to F_i(*)$ over $x_i$ and $\sim$ is the equivalence relation of concordance.
The second isomorphism follows from \cref{lem:pi0-sphere}.
On the other hand, appealing again to \cref{lem:pi0-sphere} but now for the $\infty$-sheaf $\prod_i F_i$, we have that
$$\pi_n \Bigl(\Bigl| \prod_{i} F_i(\A^\bullet) \Bigr|, x\Bigr) \cong \pi_0 \Bigl(\prod_i F_i (S^n)_{x_i}\Bigr) / {\sim}  \cong  \Bigl(\prod_i \pi_0 F_i (S^n)_{x_i}\Bigr) / {\sim}$$
where in this case $\sim$ is the equivalence relation of component-wise concordance.
The canonical map
$$\prod_{i} \pi_0 F_i(S^n)_{x_i} / {\sim} \to  \Bigl(\prod_i \pi_0 F_i (S^n)_{x_i} \Bigr) / {\sim}$$
is indeed a bijection since, in the category of Sets, infinite products commute with taking quotients by equivalence relations.
\end{proof}

\section{What does the classifying space of an $\infty$-category classify?}\label{sec:class}
In this section, we suggest an answer to the question in the title.
This expands on earlier questions and earlier answers in Moerdijk~\cite{Moerdijk} and Weiss~\cite{Weiss.Classifying}.
Even earlier results on concordance classification
of $\sC$-bundles on manifolds (or paracompact spaces) for a topological category $\sC$ (or even just a simplicial space) are due to Segal~\cite{Segal} and Stasheff~\cite{Stasheff}.
Our point of view is particularly close to Segal's.

Throughout this section, we will take $\sC$ to be a Segal space, although the discussion holds more generally for any simplicial space.
For convenience, we assume that $\sC$ is Reedy fibrant as a simplicial space, otherwise the mapping spaces below need to be derived.
(For definitions and more explanations, see Rezk~\cite{Rezk.Model}.)
For example, $\sC$ could be the (Reedy fibrant replacement of the) nerve of a (topological) category.
Informally, the following data should produce something deserving the name of a \emph{$\sC$-bundle on a manifold} $M$:
\begin{itemize}
\item an open cover $\cU = \{U_i\}$ of $M$ and a total order on its indexing set $I$,
\item maps $\{\phi_i : U_i \to \sC_0=\mathbf{ob}(\sC)\}$,
\item maps $\{ \phi_{i<j} : U_i \cap U_j \to \sC_1=\mathbf{mor}(\sC)\}$,
\item etc.
\end{itemize}
These data are then required to satisfy compatibility conditions; e.g., for a point $x \in U_i \cap U_j$, $\phi_{i < j}(x)$ is a morphism in $\sC$ from $\phi_i(x)$ to $\phi_j(x)$.
As everywhere else in this paper, \emph{space} means \emph{simplicial set}, so in the above a map from $U_i$ is taken to mean a map of simplicial sets from the singular simplices of $U_i$ to a given simplicial set.

\medskip
We make the above informal description a $\sC$-bundle precise as follows.
\begin{definition}
A \emph{$\sC$-bundle} is an open cover $\cU = \{U_i \to M\}_{i \in I}$
(we stress that here we do not require that $I$ be totally ordered)
with a simplicial space map $N \cU \to \sC$, where $N\cU$ denotes the nerve of the following topological poset.
The space of objects is
$$\bigsqcup_{\varnothing \neq S \subset I} U_S$$
where the coproduct runs over nonempty finite subsets $S$ of $I$ and $U_S \coloneq \cap_{s \in S} U_s$.
Given objects $(R, x)$ and $(S, y)$, with $x \in U_R$ and $y \in U_S$, there is a morphism $(R, x) \to (S,y)$ if and only if $R \subset S$ (so that $U_S \subset U_R$) and $x = y$.
Therefore, the space of morphisms is
$$\bigsqcup_{\varnothing \neq R \subset S} U_{S}.$$
We view $N \cU$ as a simplicial space.
Since $N \cU$ is Reedy cofibrant, the mapping space $\map(N\cU, \sC)$ agrees with the derived mapping space $\Rmap(N\cU, \sC)$.
\end{definition}

\begin{remark}
The informal description can be viewed as a special case of the definition by setting the images of certain morphisms---prescribed according to the ordering of $I$---to be identities.
Conversely, given a $\sC$-bundle it is sometimes possible to construct a $\sC$-bundle as in the informal description above
by adding to the collection $\cU$ all finite intersections of open sets in the original cover and choosing a total ordering on the resulting collection.
\end{remark}

\def\Cov{{\sf Cov}}
\def\Covsd{\Cov_{\sf sd}}

\medskip
We now build a space of $\sC$-bundles.
First, a definition:

\begin{definition}
For a manifold $M$, we define a simplicially enriched category $\Cov(M)$ of open covers $\cU$ of $M$ and their refinements.
Recall that a refinement $\cU \to \cV$ is a choice of function $\alpha : I \to J$ between the indexing sets of the covers such that $U_i \subset V_{\alpha(i)}$ for each $i \in I$.
We define a $k$-simplex in the space of morphisms of $\Cov(M)$,
$$\map(\cU, \cV)$$
to be a $(k+1)$-tuple of refinements $\alpha_0, \ldots, \alpha_k : \cU \to \cV$.
The face and degeneracy maps are clear.

The space $\map(\cU, \cV)$ may of course be empty.
If it is nonempty, it is the nerve of a groupoid, and for every pair of objects $\alpha_0, \alpha_1$ there is by construction a unique morphism $\alpha_0 \to \alpha_1$.
It follows that every $k$-sphere in $\map(\cU, \cV)$ has an unique filler, for every $k \geq 0$.
Therefore, $\map(\cU, \cV)$ is either empty or contractible.
As such, $\Cov(M)$ is equivalent (as a simplicially enriched category) to the preorder of open covers $\cU$ of $M$ with order relation $\cU \leq \cV$ if $\cU$ refines $\cV$.
\end{definition}

The assignment $\cU \mapsto N \cU$ defines a simplicially enriched functor from $\Cov(M)$ to the category of simplicial spaces, since
$\map(\cU, \cV)$ is a subspace of the space $\map(N \cU , N \cV)$ of simplicial space maps.
Indeed, each refinement $\cU \to \cV$ defines a map of simplicial spaces $N \cU \to N \cV$.
For each pair of refinements $\alpha_0, \alpha_1 : \cU \to \cV$,
the relations $U_i \subset V_{\alpha_0(i)}$ and $U_i \subset V_{\alpha_1(i)}$ imply that $U_i \subset V_{\alpha_0(i)} \cap V_{\alpha_1(i)}$
and, as such, define a simplicial map $N \cU \times \Delta[1] \to N \cV$.
More generally, a choice of refinements $\alpha_0, \ldots, \alpha_k : \cU \to \cV$ implies the relation $U_i \subset V_{\alpha_0(i)} \cap \cdots \cap V_{\alpha_k(i)}$ and so defines a map $N \cU \times \Delta[k] \to N \cV$.

\begin{definition}\label{defn:sheaf-bundles}
The $\infty$-sheaf of $\sC$-bundles is the functor which to a manifold $M$ associates the space
$$\sC(M) \coloneq \hocolim_{\cU \to M} \map(N\cU, \sC)$$
given by the homotopy colimit of the enriched functor $\cU \mapsto \map(N\cU, \sC)$ on $\Cov(M)$.
This functor is indeed enriched since on morphisms it is the restriction of the canonical map of spaces
$$\map(N \cU, N \cV) \to \map_{\Sp}(\map(N \cV, \sC), \map(N \cU, \sC))$$
to $\map_{\Cov(M)}(\cU, \cV)$.
\end{definition}

The formula in this definition applies even if $M$ has corners.
So we may view $\sC$ as a functor on the larger category of manifolds with corners and smooth maps.
In this setting, the subsheaf of sets $\bDelta^n \subset \A^n$ of \cref{sec:closedsimplices} is representable.

\begin{proposition}\label{prop:refine}
There is a canonical weak equivalence of simplicial spaces $\Ex^{\infty}\sC \to \sC(\bDelta^\bullet)$.
\end{proposition}

\begin{proof}
Let $\Cov$ denote the category $\Cov(\bDelta^n)$ of open covers of $\bDelta^n$ and refinements.
Let $\Covsd$ be the full subcategory of $\Cov$ spanned by open covers by open stars of the vertices of some barycentric subdivision of $\bDelta^n \subset \A^n$.
The set of objects of $\Covsd$ is therefore identified with the nonnegative integers: for each $i \geq 0$, the corresponding open cover $\cU(i)$ of $\bDelta^n$ is indexed by the set of vertices of the $i$-${th}$ barycentric subdivision of $\bDelta^n$.
The simplicial space $N\cU(i)$ is degreewise weakly equivalent to the simplicial discrete space $\Sd^{i+1} \Delta[n]$.
To see this, note that for a subset $S \subset \Sd^i \Delta[n]_0$ the space $U_S$ is the open star of the unique nondegenerate simplex in $\Sd^i \Delta[n]$ with vertex set $S$,
if that simplex exists, and otherwise is empty; and the $0$-simplices of $\Sd^{i+1} \Delta[n]$ are by definition the nondegenerate simplices of $\Sd^i \Delta[n]$.

For each $i \geq 0$, there is a contractible choice of morphisms $\cU({i+1}) \to \cU({i})$ in $\Covsd$.
Among these, we are interested in a specific morphism, namely the one whose underlying function between indexing sets $
\Sd^{i+1} \Delta[n]_0 = \Sd (\Sd^{i} \Delta[n])_0 \to \Sd^i \Delta[n]_0$
is the last vertex map.
The corresponding functor $\NN \to \Covsd$ that selects these morphisms is an equivalence of simplicial categories.

Write $j : \Covsd \hookrightarrow \Cov$ for the inclusion.
Clearly, every open cover of $\bDelta^n$ can be refined by one in $\Covsd$.
That is to say, for every open cover $\cV$ of $\bDelta^n$, the comma category $j/\cV$ is nonempty.
The category $j/\cV$ is equivalent to the discrete category (preorder) of open covers $\cU(i)$ in $\Covsd$ such that $\cU(i) \leq \cV$ with refinement relation $\leq$.
Clearly, $\cU(i) \leq \cU({i'})$ if and only if $i \geq i'$.
From this description it is clear that $j/\cV$ is contractible.
This shows that $j$ is homotopy final, i.e., that the homotopy colimit defining $\sC(\bDelta^n)$ may be indexed by the smaller $\Covsd$.

To summarize, we have weak equivalences
$$\hocolim_{i > 0} \map(\Sd^i \Delta[n], \sC) \to \hocolim_{\cU \in \Covsd} \map(N\cU, \sC) \xrightarrow{j^*} \sC(\bDelta^n)$$
which are functorial in $n$, and so the result follows.
\end{proof}

\begin{theorem}
For every smooth manifold $M$, the natural map
$$\sB\sC(M) \to \Rmap(M, \tB\sC)$$
is a weak equivalence.
Here $\tB\sC$ denotes the classifying space of $\sC$, i.e., the geometric realization of $\sC$ viewed as a simplicial space, and $\sB\sC$ is the functor $\sB$ applied to the $\infty$-sheaf in \cref{defn:sheaf-bundles}.
\end{theorem}
\begin{proof}
This is immediate from \cref{mainthm}, together with the identification of $|\sC(\bDelta^\bullet)|$
with $| \Ex^{\infty} \sC |$ from \cref{prop:refine}, and $| \Ex^{\infty} \sC |$ with $\tB \sC = | \sC |$ from \cref{prop:Exinfty}.
\end{proof}

\begin{remark}
Clearly, if $\sC \to \sD$ is a map inducing a weak equivalence between classifying spaces $\tB\sC \to \tB\sD$, then $\sB\sC(M) \to \sB\sD(M)$ is a weak equivalence for every manifold $M$.
This is the case, for example, if $\sC \to \sD$ is a Dwyer--Kan equivalence of Segal spaces.
\end{remark}

\appendix

\section{Technical lemmas on simplicial sets and spaces}\label{sec:appendix}

This appendix contains characterizations of weak equivalences in simplicial sets and simplicial spaces in the form that is needed in the paper.
For simplicial sets, this is classical but we could not find the statements that we need in the literature (\cref{expl:duggerisaksen,expl:duggerisaksen3}).
For simplicial spaces, this is less standard, although it may well be known.


\subsection{Special criteria for simplicial weak equivalences}

\begin{proposition}[{Dugger--Isaksen~\cite[Proposition 4.1]{Dugger.Isaksen}}]\label{prop:duggerisaksen}
A map $f : X \to Y$ between Kan complexes is a weak equivalence if and only if for every $n \geq 0$ and every commutative square
\begin{equation*}
\begin{tikzpicture}[descr/.style={fill=white}, baseline=(current bounding box.base)]
\matrix(m)[matrix of math nodes, row sep=2.5em, column sep=2.5em,
text height=1.5ex, text depth=0.25ex] {
\partial \Delta^n & X \\
\Delta^n & Y \\
};
\path[->,font=\scriptsize]
	(m-1-1) edge node [auto] {} (m-1-2)
	(m-2-1) edge node [auto] {} (m-2-2)
	(m-1-1) edge node [left] {} (m-2-1)
	(m-1-2) edge node [auto] {$f$} (m-2-2);
\path[->,dashed,font=\scriptsize]
	(m-2-1) edge node [auto] {} (m-1-2);
\end{tikzpicture}
\end{equation*}
there is a lift as pictured making the upper triangle commute and the lower triangle commute up to a homotopy $H : \Delta^n \times \Delta^1 \to Y$ which is fixed on $\partial \Delta^n$.
\end{proposition}
In Dugger--Isaksen's terminology, a map $f$ solving the lifting problem of this proposition is said to have the \emph{relative homotopy lifting property} (RHLP) with respect to $\partial \Delta^n \to \Delta^n$.
\medskip

It will be useful to think of these lifting properties in the following way.
Let $\Sp^{[1]}$ denote the category whose objects are maps of simplicial sets and morphisms are commutative squares.
Let $\tau$ denote the morphism in $\Sp^{[1]}$
\begin{equation*}
\begin{tikzpicture}[descr/.style={fill=white}, baseline=(current bounding box.base)]
\matrix(m)[matrix of math nodes, row sep=2.5em, column sep=2.5em,
text height=1.5ex, text depth=0.25ex] {
\partial \Delta^n & \Delta^n \\
\Delta^n & \Delta^n \times \Delta^1 \sqcup_{\partial \Delta^n \times \Delta^1} \partial \Delta^n \ \\
};
\path[->,font=\scriptsize]
	(m-1-1) edge node [auto] {} (m-1-2)
	(m-2-1) edge node [auto] {} (m-2-2)
	(m-1-1) edge node [left] {$i$} (m-2-1)
	(m-1-2) edge node [auto] {$j$} (m-2-2);
\end{tikzpicture}
\end{equation*}
(with source $i$ and target $j$).
\cref{prop:duggerisaksen} then reads: a map $f$ between Kan complexes is a weak equivalence if and only if $\tau^* : \map(j, f) \to \map(i, f)$ is a surjection.

\begin{corollary}\label{cor:duggerisaksen}
Let $\tau' : i' \to j'$ be a commutative square weakly equivalent to $\tau$.
(That is, $\tau'$ is related to $\tau$ by a zigzag of weak equivalences of squares.)
A map $f : X \to Y$ of simplicial sets is a weak equivalence if and only if
$$\Rmap(j', f) \to \Rmap(i', f)$$
is a surjection on $\pi_0$, where $\Rmap$ refers to the homotopy function complex in $\Sp^{[1]}$ relative to objectwise weak equivalences.
\end{corollary}
\begin{proof}
Since derived mapping spaces are invariant by weak equivalences by definition or construction, it suffices to prove that $f$ is a weak equivalence if and only if
\begin{equation}\label{eq:Rtau}
\Rmap(j, f) \to \Rmap(i, f)
\end{equation}
is a surjection on $\pi_0$.
To interpret the derived mapping spaces, let us equip $\Sp^{[1]}$ with the \emph{projective} model structure.
In this model structure, an object (i.e., map) is fibrant if source and target are Kan simplicial sets.
Without loss of generality, we may assume that $f$ is fibrant.
Cofibrant objects are simplicial maps that are cofibrations (between cofibrant objects, which is no condition here).
Cofibrations are commutative squares
\begin{equation*}
\begin{tikzpicture}[descr/.style={fill=white}, baseline=(current bounding box.base)]
\matrix(m)[matrix of math nodes, row sep=2.5em, column sep=2.5em,
text height=1.5ex, text depth=0.25ex] {
A & A' \\
B & B' \ \\
};
\path[->,font=\scriptsize]
	(m-1-1) edge node [auto] {} (m-1-2)
	(m-2-1) edge node [auto] {} (m-2-2)
	(m-1-1) edge node [left] {$i$} (m-2-1)
	(m-1-2) edge node [auto] {$j$} (m-2-2);
\end{tikzpicture}
\end{equation*}
(with source $i$ and target $j$) where the top map and the map
$$A' \sqcup_{A} B \to B'$$
are cofibrations of simplicial sets.
It is not difficult to see that the morphism $\tau$ is then a cofibration between cofibrant objects.
It follows that
$$\tau^* : \map(j, f) \to \map(i, f)$$
is identified with (\ref{eq:Rtau}) and is a Kan fibration.
Since a Kan fibration is surjective if and only if it is surjective on $\pi_0$, the result follows.
\end{proof}

Below are three examples which give rise to equivalent lifting problems:
\begin{example}\label{expl:duggerisaksen2}
Let $\tau'$ be the morphism in $\Sp^{[1]}$
\begin{equation*}
\begin{tikzpicture}[descr/.style={fill=white}, baseline=(current bounding box.base)]
\matrix(m)[matrix of math nodes, row sep=2.5em, column sep=2.5em,
text height=1.5ex, text depth=0.25ex] {
\partial \Delta^n & \partial \Delta^n \times \Delta^1 \sqcup_{\partial \Delta^n \times \{1\}} \Delta^n\\
\Delta^n & \Delta^n \times \Delta^1.\\
};
\path[->,font=\scriptsize]
	(m-1-1) edge node [auto] {} (m-1-2)
	(m-2-1) edge node [auto] {} (m-2-2)
	(m-1-1) edge node [left] {$i'$} (m-2-1)
	(m-1-2) edge node [auto] {$j'$} (m-2-2);
\end{tikzpicture}
\end{equation*}
Then $\tau'$ is weakly equivalent to $\tau$ and is a projective cofibration.
\end{example}

\begin{example}\label{expl:duggerisaksen}
Let $\tau'$ be the morphism in $\Sp^{[1]}$
\begin{equation*}
\begin{tikzpicture}[descr/.style={fill=white}, baseline=(current bounding box.base)]
\matrix(m)[matrix of math nodes, row sep=2.5em, column sep=2.5em,
text height=1.5ex, text depth=0.25ex] {
\partial \Delta^n & \partial \Delta^n \times \Delta^1 \sqcup_{\partial \Delta^n \times \{1\}} \Delta^n \\
\Lambda^{n+1} & \Lambda^{n+1} \times \Delta^1 \sqcup_{\Lambda^{n+1}  \times \{1\}} \Delta^{n+1}. \\
};
\path[->,font=\scriptsize]
	(m-1-1) edge node [auto] {} (m-1-2)
	(m-2-1) edge node [auto] {} (m-2-2)
	(m-1-1) edge node [left] {$i'$} (m-2-1)
	(m-1-2) edge node [auto] {$j'$} (m-2-2);
\end{tikzpicture}
\end{equation*}
Then $\tau'$ is weakly equivalent to $\tau$ and is a projective cofibration.
\end{example}

\begin{example}\label{expl:duggerisaksen3}
Let $D$ be the simplicial set defined as the quotient $\Delta^2 / d_0$ where $d_0 : \Delta^1 \to \Delta^2$ is the face that misses $0$.
The two remaining faces $d_1, d_2$ give two inclusions $\Delta^1 \to D$.
Let $\tau'$ be the morphism in $\Sp^{[1]}$
\begin{equation*}
\begin{tikzpicture}[descr/.style={fill=white}, baseline=(current bounding box.base)]
\matrix(m)[matrix of math nodes, row sep=2.5em, column sep=2.5em,
text height=1.5ex, text depth=0.25ex] {
\partial \Delta^n & \Delta^n \sqcup_{\partial \Delta^n} \partial \Delta^n \times \Delta^1 \\
\Delta^{n} & \Delta^{n} \times \Delta^1 \sqcup_{\partial \Delta^n \times \Delta^1} \partial \Delta^n \times D.\\
};
\path[->,font=\scriptsize]
	(m-1-1) edge node [auto] {} (m-1-2)
	(m-2-1) edge node [auto] {} (m-2-2)
	(m-1-1) edge node [left] {$i'$} (m-2-1)
	(m-1-2) edge node [auto] {$j'$} (m-2-2);
\end{tikzpicture}
\end{equation*}
Then $\tau'$ is weakly equivalent to $\tau$ and is a projective cofibration.
\end{example}

So, in view of the previous result, a map $f : X \to Y$ between Kan complexes is a weak equivalence if and only if
$$(\tau')^* : \map(j', f) \to \map(i', f)$$
is surjective for $\tau' : i' \to j'$ as in the examples above.

\subsection{Criteria for realization weak equivalences}
\begin{definition}
A simplicial space is a contravariant functor from $\Delta$ to spaces (alias simplicial sets).
\end{definition}

A simplicial space $[m] \mapsto X_m$ may be viewed as a bisimplicial set, i.e., a contravariant functor from $\Delta \times \Delta$ to $\Sets$.
However, the two $\Delta$ directions play different roles and it is important to distinguish them.

A map $X \to Y$ between simplicial spaces is a (degreewise) \emph{weak equivalence} if for each $m \geq 0$ the map $X_m \to Y_m$ is a weak equivalence of spaces.
We write $\Rmap(X,Y)$ for the homotopy function complex with respect to degreewise weak equivalences.
This may be computed as $\map(X^c, Y^f)$ in a model structure on simplicial spaces with levelwise weak equivalences, for a cofibrant replacement $X^c \to X$ and a fibrant replacement $Y \to Y^f$.
There are two canonical choices for such a model structure: the Reedy (= injective) model structure and the projective model structure.

The diagonal functor $d : \sSp \to \Sp$ has a left adjoint $d_!$ which is the unique colimit-preserving functor with $d_! (\Delta^n) =  \Delta^n \otimes \Delta[n]$.
(For a simplicial set $K$ and a simplicial space $X$, the tensor $K \otimes X$ is the simplicial space with $n$-simplices $K \times X_n$.)

There is another colimit-preserving functor $\delta : \Sp \to \sSp$ defined by $\delta(\Delta^n) = \Delta[n]$, i.e., pullback along the projection onto the first factor $\Delta \times \Delta \to \Delta$.
The projection $\Delta^n \otimes \Delta[n] \to \Delta^0 \otimes \Delta[n]$ induces a natural transformation $d_! \to \delta$.

\begin{lemma}\label{lem:d!}
For each simplicial set $X$, the natural map $d_!(X) \to \delta(X)$ is a degreewise weak equivalence of simplicial spaces.
\end{lemma}
\begin{proof}
For representables, this is clear.
A general simplicial set $X$ is a filtered colimit of finite-dimensional simplicial sets $X_i$ and filtered colimits of simplicial spaces are homotopy colimits, so it is enough to prove the statement for finite-dimensional simplicial sets.
Suppose that we have proved the statement for all simplicial sets of dimension $< n$.
We want to prove the statement for a simplicial set $X$ of dimension $n$.
Let $\sk_{n-1} X$ denote the $(n-1)^{th}$ skeleton of $X$, so that we have a pushout
$$\begin{tikzpicture}[descr/.style={fill=white}, baseline=(current bounding box.base)]
\matrix(m)[matrix of math nodes, row sep=2.5em, column sep=2.5em,
text height=1.5ex, text depth=0.25ex] {
\bigsqcup_{X_n} \partial \Delta^n  & 	\sk_{n-1} X  \\
\bigsqcup_{X_n} \Delta^n & X  \\
};
\path[->,font=\scriptsize]
	(m-1-1) edge node [auto] {} (m-1-2)
	(m-2-1) edge node [auto] {} (m-2-2)
	(m-1-1) edge node [left] {} (m-2-1)
	(m-1-2) edge node [auto] {} (m-2-2);
\end{tikzpicture}$$
Since $d_!$ and $\delta$ are colimit-preserving, the result follows by induction and the case of representables.
\end{proof}

\begin{proposition}\label{prop:DuggerIsaksenSS2}
Let $f : X \to Y$ be a map between Reedy fibrant simplicial spaces which satisfy the Kan condition.
Then $|f| : |X| \to |Y|$ is a weak equivalence if and only if every square
$$\begin{tikzpicture}[descr/.style={fill=white}, baseline=(current bounding box.base)]
\matrix(m)[matrix of math nodes, row sep=2.5em, column sep=2.5em,
text height=1.5ex, text depth=0.25ex] {
\partial \Delta[n] & X \\
\Delta[n] & Y \\
};
\path[->,font=\scriptsize]
	(m-1-1) edge node [auto] {} (m-1-2)
	(m-2-1) edge node [auto] {} (m-2-2)
	(m-1-1) edge node [left] {} (m-2-1)
	(m-1-2) edge node [auto] {} (m-2-2);
\path[dashed,->,font=\scriptsize]
(m-2-1) edge node [left] {} (m-1-2);
\end{tikzpicture}$$
has a lift as pictured making the lower triangle commute up to a given homotopy $\Delta[n] \times \Delta[1] \to Y$ and making the upper triangle commute up to a given homotopy $\partial \Delta[n] \times \Delta[1] \to X$.
These two homotopies are required to be homotopic as maps $\partial \Delta[n] \times \Delta[1] \to Y$ and the homotopy should be constant on $\partial \Delta[n] \times \partial \Delta[1]$.
\end{proposition}

\begin{proof}
Since $X$ and $Y$ are Kan complexes and $d$ preserves Kan fibrations (\cref{rem:strongKan}), $|X|$ and $|Y|$ are Kan complexes.
In view of \cref{prop:duggerisaksen}, \cref{expl:duggerisaksen3} and the remarks that follow it, $|f|$ is a weak equivalence if and only if the map
$${\tau'}^* : \map(i', |f|) \to \map(j', |f|)$$
is surjective (using the notation from \cref{expl:duggerisaksen3}).
By adjunction, this is equivalent to saying that
\begin{equation}\label{eq:dmap1}
{(d_! \tau')}^* : \map(d_!i', f) \to \map(d_! j', f)
\end{equation}
is surjective.
Since $f$ is a map between Reedy fibrant simplicial spaces,
it is a fibrant object in $\sSp^{[1]}$ with the \emph{projective} model structure on the category of functors $[1] \to \sSp$,
where $\sSp$ is equipped with the Reedy model structure.
Since $d_!$ preserves monomorphisms, $d_! \tau$ is a cofibration between cofibrant objects in that same model structure
(see the proof of \cref{cor:duggerisaksen}).
Therefore, the map (\ref{eq:dmap1}) is a fibration and as such it is surjective if and only if it is surjective on $\pi_0$.
These considerations also lead us to identify (\ref{eq:dmap1}) with the map on derived mapping spaces
$${(d_! \tau')}^* : \Rmap(d_!i', f) \to \Rmap(d_! j', f),$$
which by \cref{lem:d!} is identified with
$$\delta(\tau')^* : \Rmap(\delta(i'), f) \to \Rmap(\delta(j'), f).$$
Since $\delta(\tau')$ is also a cofibration between cofibrant objects, this map is identified with
$$\delta(\tau')^* : \map(\delta(i'), f) \to \map(\delta(j'), f).$$
The surjectivity of this last map is equivalent to the existence of the lift in the statement of the proposition.
\end{proof}

\begin{corollary}\label{cor:DuggerIsaksenSS2}
Let $f : X \to Y$ be a map between Reedy fibrant simplicial spaces.
Suppose that for every $j \geq 0$ and every square
$$\begin{tikzpicture}[descr/.style={fill=white}, baseline=(current bounding box.base)]
\matrix(m)[matrix of math nodes, row sep=2.5em, column sep=2.5em,
text height=1.5ex, text depth=0.25ex] {
\Sd^j \partial \Delta[n] & X \\
\Sd^j \Delta[n] & Y \\
};
\path[->,font=\scriptsize]
	(m-1-1) edge node [auto] {} (m-1-2)
	(m-2-1) edge node [auto] {} (m-2-2)
	(m-1-1) edge node [left] {} (m-2-1)
	(m-1-2) edge node [auto] {} (m-2-2);
\path[dashed,->,font=\scriptsize] (m-2-1) edge node [left] {} (m-1-2);
\end{tikzpicture}$$
there is a lift as pictured making the lower triangle commute up to a given homotopy $\Sd^j \Delta[n] \times \Delta[1] \to Y$ and making the upper triangle commute up to a given homotopy $\Sd^j \partial \Delta[n] \times \Delta[1] \to X$.
These two homotopies are required to be homotopic as maps $\Sd^j \partial \Delta[n] \times \Delta[1] \to Y$ and the homotopy should be constant on $\Sd^j \partial \Delta[n] \times \partial \Delta[1]$.
Then $|f| : |X| \to |Y|$ is a weak equivalence.
\end{corollary}

\begin{proof}
Apply \cref{prop:DuggerIsaksenSS2} replacing $X$ and $Y$ by the simplicial spaces $\Ex^\infty X$ and $\Ex^\infty Y$, which satisfy the Kan condition by \cref{prop:Exinfty}.
\end{proof}

\subsection{Properties of subdivisions of simplicial spaces}\label{sec:appendixex}

Recall the simplicial subdivision $\sd \Delta^n$, i.e., the nerve of the poset of nonempty subsets of $[n]= \{0, \ldots, n\}$.

\begin{definition}
\label{defn:ex}
Denote by $$\Sd:\sSp\to\sSp$$
the simplicial left adjoint functor
that sends $\Delta[n]$ to $\sd \Delta^n$ viewed as a simplicial \emph{discrete} space.
Denote by $$\Ex:\sSp\to\sSp$$
the simplicial right adjoint functor of $\Sd$.
\end{definition}

Every simplicial discrete space is Reedy cofibrant, so replacing $X$ by a Reedy fibrant simplicial space $X^f$,
we may write the right derived functor of $\Ex$ evaluated at~$X$ as the (honest) mapping space $\map(\Sd\Delta[n], X^f)$.

There is a natural map $\gamma : \Sd \Delta[n] \to \Delta[n]$, sending a subset $\{i_0, \ldots, i_k\} \subset [n]$ to $i_k$ (the last vertex).
The colimit
$$X \xrightarrow{\gamma^*} \Ex X \xrightarrow{\gamma^*} \Ex^2 X \to \cdots$$
is denoted by $\Ex^\infty X$.
The map $\gamma$ has a section $\Delta[n] \to \Sd \Delta[n]$ so if $X$ is Reedy fibrant, all the maps in the tower are degreewise cofibrations and so the colimit computes the homotopy colimit.

We collect the important properties of the $\Ex^\infty$ endofunctor below.
These parallel (or, rather, include) the well-known ones for simplicial \emph{sets}.
\begin{proposition}\label{prop:Exinfty}
For a simplicial space $X$:
\begin{enumerate}
\item $\Ex^\infty X$ is a Kan simplicial space.
\item $X \to \Ex^\infty X$ is a weak equivalence after geometric realization.
\item For each $i$, including $i = \infty$, $\Ex^i$ preserves (trivial) Kan fibrations, zero simplices and finite homotopy limits.
\end{enumerate}
\end{proposition}
\begin{proof}
By construction, the functor $\Ex^{i}$, for $0 \leq i \leq \infty$, sends weak equivalences of simplicial spaces to weak equivalences.
If $X$ is Reedy fibrant then
$$\map(\Sd \Delta[n], X) \to \map(\Sd \partial \Delta[n], X)$$
is a fibration (since $\Sd \partial \Delta[n] \to \Sd \Delta[n]$ is a degreewise monomorphism, hence a cofibration).
Therefore, $\Ex X$ is Reedy fibrant.
By standard compactness arguments, it follows that $\Ex^\infty X$ is also Reedy fibrant.
Hence, in proving $(1)$, $(2)$ and $(3)$, we may assume from the outset that $X$ is Reedy fibrant.

The arguments to prove $(1)$ and $(3)$ are identical to the classical ones for simplicial sets, so we do not reproduce them here.
As for $(2)$, take a trivial Kan fibration $X' \to X$ where $X'$ is a simplicial \emph{set} (see Lurie~\cite[Proposition 7]{Lurie.Simplicial}), and consider
the square
\begin{equation*}
\begin{tikzpicture}[descr/.style={fill=white}, baseline=(current bounding box.base)]
\matrix(m)[matrix of math nodes, row sep=2.5em, column sep=2.5em,
text height=1.5ex, text depth=0.25ex] {
X' & \Ex X' \\
X & \Ex X. \\
};
\path[->,font=\scriptsize]
	(m-1-1) edge node [auto] {} (m-1-2)
	(m-2-1) edge node [auto] {} (m-2-2)
	(m-1-1) edge node [left] {} (m-2-1)
	(m-1-2) edge node [auto] {} (m-2-2);
\end{tikzpicture}
\end{equation*}
Since the diagonal preserves trivial Kan fibrations, the vertical maps are weak equivalences after applying the diagonal (for the right-hand one, use part $(3)$).
The top horizontal map is a weak equivalence, e.g., see Goerss--Jardine~\cite[III.4.6]{Goerss.Jardinee}.
We conclude that the diagonal of the lower map is a weak equivalence.
\end{proof}



\end{document}